\pdfoutput=1
\documentclass{amsart}

\usepackage{xcolor}

\usepackage{amssymb}
\usepackage{amsthm}
\usepackage{delimset}
\usepackage{enumitem}
\usepackage{thmtools}
\usepackage{thm-restate}
\usepackage{hyperref}
\usepackage{cleveref}

\setlist[enumerate,1]{label={\upshape(\roman*)}}
\setlist[enumerate,2]{label={\upshape(\alph*)}}

\theoremstyle{plain}
\newtheorem{theorem}{Theorem}[section]
\newtheorem{lemma}[theorem]{Lemma}
\newtheorem{proposition}[theorem]{Proposition}

\theoremstyle{definition}
\newtheorem{definition}[theorem]{Definition}
\newtheorem{remark}[theorem]{Remark}

\newtheorem{df}[theorem]{Definition}
\newtheorem{ex}[theorem]{Example}

\newcommand{\term}[1]{\emph{#1}} % for defining terms
\newcommand{\defeq}{=_{df}} % definitional equality
\newcommand{\A}{\mathcal{A}}
\newcommand{\B}{\mathcal{B}}
\newcommand{\C}{\mathcal{C}}
\newcommand{\D}{\mathcal{D}}
\newcommand{\E}{\mathcal{E}}
\newcommand{\F}{\mathcal{F}}
\newcommand{\J}{\mathcal{J}}
\newcommand{\W}{\mathcal{W}}

\newcommand{\K}{\mathcal{K}} % compact operators \ell^2

\newcommand{\Mul}{\mathcal{M}} % Multiplier algebra
\newcommand{\paschkedual}[2]{#1^d_{#2}}

\numberwithin{equation}{section}

\DeclareMathOperator{\diam}{diam}
\DeclareMathOperator{\Ped}{Ped}
\DeclareMathOperator{\Ideal}{Ideal}

\DeclareMathOperator{\supp}{supp}
\DeclareMathOperator{\osupp}{osupp}
\DeclareMathOperator{\Ext}{Ext}

\DeclareMathOperator{\trace}{tr}
\DeclareMathOperator{\ind}{ind}
\DeclareMathOperator{\id}{id}
\DeclareMathOperator{\Aut}{Aut} % automorphism group
\DeclareMathOperator{\spec}{sp} % spectrm

\newcommand{\cstar}{\textrm{C}*\kern-0.15ex} % C* in text, with negative kerning
\newcommand{\paue}{\approxeq} % proper asymptotic unitary equivalence
\newcommand{\nat}{\mathbb{N}} % natural numbers
\newcommand{\real}{\mathbb{R}} % real numbers
\newcommand{\complex}{\mathbb{C}} % complex numbers
\newcommand{\mat}{M} % matrix algebra
\renewcommand{\k}{K} % K-theory
\newcommand{\kk}{KK} % KK-theory
\newcommand{\kkz}{KK^0} % KK^0
\newcommand{\kko}{KK^1} % KK^1
\newcommand{\kz}{K_0} % K_0
\newcommand{\ko}{K_1} % K_1
\newcommand{\kkhig}{\kk_{\mathrm{Higson}}}
\newcommand{\ec}{\delimpair[{[.]:}]}

\title[Paschke dual algebra]{$\ko$-injectivity of the Paschke dual algebra for certain
simple C*-algebras}

\author{Jireh Loreaux}

\address{Department of Mathematics and Statistics\\
Southern Illinois University Edwardsville\\
1 Hairpin Dr.\\
Edwardsville, IL\\
62026-1653\\ USA}

\email{jloreau@siue.edu}

\author{P. W. Ng}

\address{Department of Mathematics\\
University of Louisiana at Lafayette\\
217 Maxim Doucet Hall\\
P. O. Box 43568\\
Lafayette, Louisiana\\
70504--3568\\
USA}

\email{png@louisiana.edu}

\author{Arindam Sutradhar}

\address{Department of Mathematics\\
University of Louisiana at Lafayette\\
217 Maxim Doucet Hall\\
P. O. Box 43568\\
Lafayette, Louisiana\\
70504--3568\\
USA}

\email{arindam.sutradhar1@louisiana.edu}

\begin{document}

\begin{abstract}
Let $\B$ be a nonunital separable simple stable
 C*-algebra with strict comparison of positive
elements and $T(\B)$ having finite extreme boundary, and let $\A$ be a simple unital separable 
nuclear C*-algebra.  We prove that the Paschke dual algebra $\A^d_{\B}$ is $\ko$-injective.
  
As a consequence, we obtain interesting $\kk$-uniqueness theorems which generalize the Brown--Douglas--Fillmore essential codimension property.
\end{abstract}

\maketitle

\section{Introduction}

In their seminal paper \cite[Remark~4.9]{BDF-1973-PoaCoOT}, Brown, Douglas and Fillmore
(BDF)  classified essentially normal operators using Fredholm indices.
In the course of proving functorial properties of their homology $\operatorname{Ext}(X)$, BDF introduced
the \term{essential codimension} for a pair of projections $P,Q \in B(\ell_2)$ whose
difference is compact.

Since that time, the concept of essential codimension has had numerous applications including the computation of spectral flow in semifinite von Neumann algebras \cite{BCP+-2006-Agatoeo}, as well as the explanation \cite{KL-2017-IEOT,Lor-2019-JOT} of unexpected integers arising in strong sums of projections \cite{KNZ-2009-JFA}, diagonals of projections \cite{Kad-2002-PNASUa}, and diagonals of normal operators with finite spectrum \cite{Arv-2007-PNASU,Jas-2013-JFA,BJ-2015-TAMS}.

We here present the definition of essential codimension in $B(\ell_2)$, but the translation to arbitrary semifinite factors is straightforward.

\begin{definition}
  \label{def:standard-essential-codimension}
  Given projections $P,Q \in B(\ell_2)$ with $P - Q \in \K$, the \term{essential codimension of $Q$ in $P$} is the following quantity:
  \begin{align*}
    \ec{P}{Q} &\defeq
                \begin{cases}
                  \trace(P) - \trace(Q) & \text{if}\ \trace(P) + \trace(Q) < \infty, \\[0.5em]
                  \ind(V^*W) &
                  \parbox[t]{28ex}{%
                    if $\trace(P) = \trace(Q) = \infty,$ \\
                    where $V^*V = W^*W =1,$ \\
                    $VV^* = Q, WW^* = P,$%
                  }
                \end{cases}
  \end{align*}
  where $\ind$ denotes the Fredholm index.
  We note that when $Q \le P$, $\ec{P}{Q}$ coincides with the usual codimension of $Q$ in $P$.
  Note also that $\ec{P}{Q} \in \kz(\K) \cong \mathbb{Z}$ where the isomorphism is induced by the trace.
\end{definition}

The fundamental property introduced by
Brown--Douglas--Fillmore is encapsulated in the following theorem, whose proof extends easily to semifinite factors (see \cite{BL-2012-CJM} or \cite{KL-2017-IEOT}).

\begin{theorem}
  \label{thm:bdf-unitary-equivalence-semifinite-factors}
  If $P,Q \in B(\ell_2)$ are projections with $P-Q \in \K$, then $\ec{P}{Q} = 0$ if and only if there is a unitary $U \in 1 + \K$ conjugating $P$ to $Q$, i.e., $UPU^{*} = Q$.
\end{theorem}

It has been recognized that essential codimension can
be realized as an element of $\kkz(\complex,\K)$ (or $\kkz(\complex,\complex)$) where one identifies $\ec{P}{Q}$ with the equivalence class $[\phi,\psi]$ where $\phi,\psi : \complex \to B(\ell_2)$ are *-homomorphisms with $\phi(1) = P$ and $\psi(1) = Q$.
This leads to a natural generalization of essential codimension to the setting $\kkz(\A,\B)$ where $\A$ is a separable nuclear \cstar-algebra and $\B$ is separable stable \cstar-algebra,
and thus, uniqueness results which generalize
\Cref{thm:bdf-unitary-equivalence-semifinite-factors}
(e.g., see \cite{Lee-2011-JFA}; see also, \cite{BL-2012-CJM},
\cite{Lin-2002-JOT}, \cite{DE-2001-KT},
\cite{LN-2020-IEOT}.)

It turns out that a sufficient condition for such generalizations
is the $\ko$-injectivity of the Paschke dual algebra $\paschkedual{\A}{\B}$ (see \Cref{thm:MainUniqueness}), which is a key ingredient in Paschke duality
$\kk^j(\A, \B) \cong \k_{j+1}(\paschkedual{\A}{\B})$ for $j = 0,1$
(e.g., see \cite{Pas-1981-PJM}, \cite{Tho-2001-PAMS},
\cite{Val-1983-PJM}).
Recall that a unital \cstar-algebra $\C$ is \term{$\ko$-injective} if the 
standard map from $\mathcal{U}(\C)/\mathcal{U}_0(\C)$ to $K_1(\C)$ is injective.
The Paschke dual algebra $\paschkedual{\A}{\B}$ is the relative commutant in the corona algebra $\C(\B)$ of the image of $\A$ under a strongly
unital trivial absorbing extension (see \Cref{def:paschke-dual}).

$\ko$-injectivity of the Paschke dual algebra is, in itself, an interesting question.  For example, in the case where $\A$ is unital, the Paschke dual algebra is properly infinite (e.g., see \Cref{lem:Jan1520216AM}), and it is an interesting open question of Blanchard,
Rohde and R\o{}rdam whether every properly infinite unital \cstar-algebra is $\ko$-injective (\cite[Question~2.9]{BRR-2008-JNG}).\footnote{Recall that a 
unital C*-algebra $\C$ is said to be \emph{properly infinite} if 
$1_{\C} \otimes 1_{\C} \preceq 1_{\C}$ in $\mathbb{M}_2(\C)$.}
Consequently, we focus attention on determining conditions on $\A$ and $\B$ which guarantee that the Paschke dual algebra is $\ko$-injective.

In this paper we continue the investigation \cite{LN-2020-IEOT} by the first and second named authors (itself in the spirit of \cite{Lee-2011-JFA,BL-2012-CJM,Lee-2013-JFA,Lee-2018-JMAA}) to obtain generalizations of \Cref{thm:bdf-unitary-equivalence-semifinite-factors} in various contexts.
One such generalization from \cite{LN-2020-IEOT} was the following result (see \Cref{sec:paschke-dual-uniq} and the beginning of \Cref{sec:paschke-dual-simple} for many definitions):

\begin{theorem}[\protect{\cite[Theorem~3.5]{LN-2020-IEOT}}]
  Let $\A$ be a unital separable simple nuclear \cstar-algebra, and $\B$ a separable simple stable \cstar-algebra with a nonzero projection, strict comparison of positive elements and for which $T(\B)$ (see the beginning of \Cref{sec:paschke-dual-simple}) has finitely many extreme points.
  Suppose that there exists a *-embedding $\A \hookrightarrow \B$.

  Let $\phi,\psi : \A \to \Mul(\B)$ be unital trivial full extensions such that $\phi(a) - \psi(a) \in \B$ for all $a \in \A$.
  Then $[\phi,\psi] = 0$ in $\kk(\A,\B)$ if and only if $\phi,\psi$ are properly asymptotically unitarily equivalent.
  \label{thm:Jan152021Restrictive}
\end{theorem}

See \Cref{def:proper-asymptotic-unitary-equivalence} for the definition of proper asymptotic unitary equivalence.
We will also recall much of the other terminology in the body of the paper.
Of course, the hypotheses here are quite restrictive.  For example, the condition
of $\B$ having a projection rules out many interesting C*-algebras like the
Razak algebra $\W$, a $\kk$-contractible stably projectionless nuclear
 C*-algebra with
unique trace, which has
a fundamental role in the Elliott classification program and has been very popular in recent
years.  
  Even worse,
the severe 
restriction that $\A$ be embeddable into $\B$ completely rules out the numerous
cases where $\A$ is infinite and $\B$ is stably finite.  We note that one (of many)
quick example of this is the case where $\A = O_{\infty}$ and 
$\B = C^*_r(\mathbb{F}_2) \otimes \K$.  In this case, $\A$ is simple purely infinite,
$\B$ is stably finite with unique trace and 
strict comparison of positive elements, and additionally,
$\B$ is not even Jiang--Su stable.  These examples are covered by the new theorem
in the present paper.                      

Thus, in this paper we remove many of the restraints of \Cref{thm:Jan152021Restrictive} by proving the following
result:  

\begin{restatable*}{theorem}{maintheorem}
  \label{thm:main-theorem}
  Let $\A$, $\B$ be separable simple C*-algebras with $\A$ unital and nuclear,
  with $\B$ stable, finite 
  and having strict comparison for positive elements, and
  with $T(\B)$ having finitely many extreme points.

  Let $\phi, \psi : \A \rightarrow \Mul(\B)$ be two unital full  
  extensions with $\phi(a) - \psi(a) \in \B$ for all $a \in \A$.

  Then $[\phi, \psi] = 0$ in $\kk(\A, \B)$ if and only if
  $\phi$ and $\psi$ are properly asymptotically unitarily equivalent.
\end{restatable*}

Historically, there are two main reasons for the interest in this result.
Firstly, results of this type led to the stable uniqueness theorems of classification theory,
which in turn led to bounded stable uniqueness and then non-stable uniqueness theorems, which
are key tools in the Elliott classification program.
Perhaps one definitive historical development was Lin's proof of the first general
stable uniqueness theorem, which eventually led to his TAF uniqueness theorem and the classification of simple amenable TAF algebras which axiomatized the work of Elliott and Gong (\cite[Theorem~4.3]{Lin-2002-JOT}, \cite[Theorems~2.1 and 2.3]{Lin-2001-CJM}, \cite{DE-2001-KT,EG-1996-AM}).
Secondly, the above result complements nicely the theory of extensions (or $\kko(\A,\B)$).
We note that extension theory is essentially an asymptotic theory of operators (e.g., the classification of essentially normal operators up to unitary equivalence modulo the compacts in \cite{BDF-1973-PoaCoOT}).
In contrast, the above result concerns the local theory of operators (e.g., essential codimension).

\Cref{thm:main-theorem} is, in itself, of great interest, but in a future paper, we will also 
work out some interesting consequences to projection lifting in
C*-algebras.

As previously mentioned, the key to establishing this result, as in \cite{LN-2020-IEOT}, is to prove that the Paschke dual algebra $\paschkedual{\A}{\B}$ is $\ko$-injective.
To this end, in \Cref{sec:paschke-dual-uniq}, we fix some notation on extension theory  and the Paschke dual algebra $\paschkedual{\A}{\B}$ (see \Cref{def:paschke-dual}).
We review and slightly improve upon the properties of the Paschke dual algebra developed in \cite{LN-2020-IEOT}, many of which are generalizations of results from Paschke's original paper which was for the case when $\B = \K$.
We conclude \Cref{sec:paschke-dual-uniq} with \Cref{prop:paschke-dual-purely-infinite}, which says that $\paschkedual{\A}{\B}$ is purely infinite when $\A$ is unital separable simple nuclear and either $\B = \K$ or $\B$ is separable stable simple purely infinite.
This section serves as a short summary of Paschke dual algebras, and their general properties and consequences, and as a clean reference for the future.

In Section \ref{sec:paschke-dual-simple}, which is the longest
section, we prove that the Paschke dual algebra $\paschkedual{\A}{\B}$ is $\ko$-injective when $\A,\B$ are separable and simple with $\A$ unital nuclear, $\B$ stable and having strict comparison, and $T(\B)$ having 
finite extreme boundary (see \Cref{lem:PIK1Injective} and \Cref{thm:k1-injectivity-simple-nuclear-strict-comparison}).  Then, combined with the results
from \Cref{sec:paschke-dual-uniq}, this leads to the required uniqueness
results.  Part of the goal of this section is to provide relatively general
strategies and arguments, and good examples of this are the arguments
of Subsection \ref{subsec:AsympLaurent}.   Firstly, such abstraction simplifies the technical computations and facilitates understanding.  Secondly, this will also
provide a strategy for and make easier future 
generalizations.

While this paper
utilizes $\kk$-theory and other machinery, it should be accessible to a student with
good knowledge of basic C*-algebra (say, at the level of \cite{Dav-1996}
and \cite{Weg-1993}).  We provide many definitions and references, and if such a reader 
were willing to go to the references at the appropriate moments, then he or she should be able to
work through this paper without great difficulty. In fact, while we also 
use other 
versions of $\kk$-theory in this paper, 
we focus primarily on the extension theory
and generalized homomorphism approaches.  For the student coming from analysis,
these two approaches form an excellent introduction to $\kk$-theory and 
its applications to operator theory (this is especially true of the extension
theory approach which has the closest and historically oldest connections
to operator theory; e.g., see \cite{BDF-1973-PoaCoOT}). Indeed, the subject matter of this paper
is itself a very good introduction.
%Some good references to KK theory
%are \cite{Bla-1998}, \cite{JT-1991}, \cite{Lin-2001}  and \cite{RLL-2000}.}}   

%The reader of this paper should be familiar with many aspects 
%of \cstar-algebra theory and $\k$-theory, especially extension 
%theory and $\kk$-theory.
%It would be helpful if the reader was also knowledgeable of the basic 
%ideas from the modern theory of simple nuclear \cstar-algebras.
%The reader should be comfortable with the concepts and computations in 
%the following books and the references therein:
%\cite{Dav-1996,Weg-1993,Bla-1998,JT-1991,Lin-2001,RLL-2000}.
%As examples, here are some concepts used in this paper: absorbing extensions, 
%Busby invariant, Brown--Douglas--Fillmore theory, essential codimension, 
%$\kk$-groups, $\ko$-injectivity, multiplier and corona algebras, properly 
%infinite projections, Jiang--Su stability, strict comparison, cancellation 
%of projections, nuclearity.

%Once more, good basic references that a reader should be familiar with 
%are \cite{Dav-1996} and \cite{Weg-1993}.  

Let us elaborate a bit more on the last
remarks from the previous paragraph.  
Firstly, part of the power of $\kk$-theory is that 
often the same object can be defined in many different ways, leading to many 
connections and applications.  A good
example of this is the many different versions of $\kk$-theory, which, while 
useful, 
can be bewildering to the beginning student.  
In this paper, we will be using three versions
of $\kk$-theory:  Extension theory, the generalized homomorphism version, 
and Higson's version.  
We will emphasize extension theory (though we will not define
$\Ext(\A, \B)$) which has the oldest and closest connections to
operator theory, and is a very appropriate beginning for the student of analysis. 
We focus on absorbing extensions which are related to a lot of fundamental operator
theory (e.g., Voiculescu's noncommutative Weyl--von Neumann theorem). A good basic
reference for absorbing extensions is \cite{EK-2001-PJM} -- but we also
recommend the historical works \cite{BDF-1973-PoaCoOT}, \cite{Arv-1977-DMJ}, 
\cite{Voi-1976-RRMPA}, \cite{Kas-1980-JOT} which provide a lot of insight as well
as connections to concrete operator theory.   
Nonetheless, it is important to learn the other versions of $\kk$-theory or 
at least refer to some of them when reading this paper, and for
the beginning student, we recommend \cite{Bla-1998} and \cite{JT-1991}, as well as the
original paper \cite{Kas-1981-MUI}.  [Linbook] also
provides nice information especially
 relating to uniqueness theorems.  
We also recommend \cite{Lin-2002-JOT} and \cite{Lin-2001-CJM}.

\section{Paschke duals and uniqueness}
\label{sec:paschke-dual-uniq}

We begin by fixing some notation and recalling some basic facts from extension theory.
More detailed references for extension theory can be found in \cite{Weg-1993,Bla-1998}.
Given a \cstar-algebra $\B$, we let $\Mul(\B)$ denote the multiplier algebra of $\B$, and $\C(\B) \defeq \Mul(\B)/\B$ denote the corona algebra of $\B$.
Recall that, roughly speaking, $\Mul(\B)$ is the ``largest unital C*-algebra
containing $\B$ as an essential ideal", and that $\Mul(\B)$ ``encodes" the
extension theory of $\B$.
In fact, recall that given an extension of \cstar-algebras
\begin{equation*}
  0 \to \B \to \mathcal{E} \to \A \to 0,
\end{equation*}
one associates the \term{Busby invariant}, which is a *-homomorphism $\phi : \A \to \C(\B)$.
Conversely, given such a *-homomorphism $\phi : \A \to \C(\B)$, one can obtain an extension of $\B$ by $\A$ whose Busby invariant is $\phi$.
It is well-known that the extension corresponding to a given Busby invariant is unique up to \term{strong isomorphism} (in the terminology of Blackadar; see \cite[15.1--15.4]{Bla-1998}).
Our results are all invariant under strong isomorphism, and therefore we will simply refer to $\phi : \A \to \C(\B)$ as an extension.
When $\phi$ is injective, then $\phi$ is called an \term{essential extension}.  This corresponds to $\B$ being an essential ideal of $\E$.

A \term{trivial extension} is one for which the short exact sequence is split exact, or equivalently, the Busby invariant factors through the quotient map $\pi : \Mul(\B) \to \C(\B)$ via a *-homomorphism $\phi_0 : \A \to \Mul(\B)$ so that $\phi = \pi \circ \phi_0$.
If $\phi : \A \to \C(\B)$ is a trivial extension, then by a slight abuse of terminology we also refer to the map $\phi_0$ as a trivial extension.
When $\A$ is a unital \cstar-algebra and $\phi$ is a unital *-homomorphism, we say that $\phi$ is a \term{unital extension}.
If, in addition, $\phi$ is trivial and some lift $\phi_0$ is unital, $\phi$ is said to be a \term{strongly unital trivial extension}.
In the original BDF paper [BDFPaper1], $\A = C(X)$ where $X \subset \mathbb{C}$ compact
and $\B = \K$, and the  class of trivial extensions (which is the ``index zero" case
of BDF theory) correspond to the essentially normal operators in $\mathbb{B}(l_2)$
which 
are compact perturbations of normal operators.  

Two extensions $\phi, \psi : \A \rightarrow \C(\B)$ are said to be \term{unitarily equivalent} (and denoted $\phi \sim \psi$) if there exists a unitary
$U \in \Mul(\B)$ such that
\begin{equation*}
  \pi(U) \phi(\cdot) \pi(U)^* = \psi(\cdot).
\end{equation*}
(In the literature, unitary equivalence is sometimes called \term{strong unitary equivalence}.)
In the original BDF paper, classifying extensions up to unitary equivalence corresponds to classifying essentially normal operators in $\mathbb{B}(l_2)$ up to unitary
equivalence modulo the compacts (see, for example, \cite{BDF-1973-PoaCoOT} or
\cite[Chapter IX, Section 1]{Dav-1996}).  Thus, extension theory is concerned with 
``asymptotic" features of operators.

Suppose that $\phi_0, \psi_0 : \A \rightarrow \Mul(\B)$ are *-homomorphisms.  We will sometimes write $\phi_0 \sim \psi_0$
to mean $\pi \circ \phi_0 \sim \pi \circ \psi_0$, i.e., that the corresponding extensions are unitarily equivalent.

When $\B$ is stable, there are isometries $V,W \in \Mul(\B)$ for which $VV^{*} + WW^{*} = 1$.
Given extensions $\phi,\psi : \A \to \C(\B)$, their \term{BDF sum} is the extension
\begin{equation*}
  (\phi \oplus \psi)(\cdot) \defeq \pi(V) \phi(\cdot) \pi(V^{*}) + \pi(W) \psi(\cdot) \pi(W^{*}).
\end{equation*}
It is not hard to see that the BDF sum is well-defined up to unitary equivalence (i.e., up to $\sim$).
An extension $\phi : \A \to \C(\B)$ is \term{absorbing} (respectively, \term{unitally absorbing}) if
\begin{equation*}
  \phi \sim \phi \oplus \psi
\end{equation*}
for any (respectively, \term{strongly unital}) trivial extension $\psi$.
Note that since any absorbing extension must absorb the zero *-homomorphism, it is necessarily nonunital.  The convention (which we will follow) is that if $\phi$ is a unital extension then when we say that $\phi$ is \term{absorbing}, we mean that $\phi$ is unitally absorbing.

The notion of absorbing extension originates from concrete operator theory and 
has many interesting applications (e.g., \cite{BDF-1973-PoaCoOT}, 
\cite{Arv-1977-DMJ}, \cite{Voi-1976-RRMPA}).
From a global point of view, it singles out the class of extensions that is classified by Kasparov's $\kko$ (which is a highly computable object).
In more detail,
let $\A$ be a separable nuclear C*-algebra and $\B$ be a $\sigma$-unital C*-algebra.
Kasparov proved that 
$$\kko(\A, \B) = \{ \makebox{unitary equivalence classes of
absorbing extensions } \A \rightarrow \C(\B \otimes \K) \}.$$ 
(See \cite[Theorem~1]{Kas-1981-MUI}; see also \cite{Bla-1998} 
and \cite{JT-1991}) 
Note that this also indicates that Kasparov's $\kko$ is about the ``asymptotic
aspects" of operators (e.g., classifying essentially normal operators up to
unitary equivalence modulo the compacts).  As we will see, in contrast,
Kasparov's $\kkz$ is about the ``local aspects" of operators.  This is exemplified by the generalized homomorphism picture of $\kk$-and the 
uniqueness theorem of our paper.

In \cite{Tho-2001-PAMS}, Thomsen establishes the existence of an absorbing trivial extension of $\B$ by $\A$, when $\A,\B$ are separable with $\B$ stable, and provides several equivalent characterizations of absorbing trivial extensions.
Important precursors to Thomsen's work, which are themselves of interest and give significantly more information, are contained in \cite{Voi-1976-RRMPA,Kas-1980-JOT,Lin-2002-JOT,EK-2001-PJM}.

In \cite{Pas-1981-PJM}, Paschke introduced, for a separable unital \cstar-algebra $\A$ and a unital trivial essential extension $\phi : \A \to B(\ell_2(\nat))$ (which is necessarily absorbing by Voiculescu's theorem \cite{Voi-1976-RRMPA}),
 the subalgebra $(\pi \circ \phi (\A))'$ of the Calkin algebra and proved, in the language of $\kk$-theory, $\k_j((\pi \circ \phi (\A))') \cong \kk^{j+1}(\A,\complex)$.
Soon after in \cite{Val-1983-PJM}, Valette extended these ideas and techniques to a pair of algebras $\A,\B$.
In particular, Valette proved \cite[Proposition~3]{Val-1983-PJM} that if $\A$ is a separable unital nuclear \cstar-algebra and $\B$ is a stable $\sigma$-unital \cstar-algebra, then $\k_j((\pi \circ \phi (\A))') \cong \kk^{j+1}(\A,\B)$, where $\phi : \A \to \Mul(\B)$ is a unital absorbing trivial extension.   This is the so-called \term{Paschke duality}, and it has been generalized to more general algebras (e.g., see \cite{Tho-2001-PAMS}), but in this paper we focus on the case when $\A$ is unital and nuclear.
It is for this reason that this algebra $(\pi \circ \phi (\A))'$ gets its name, the \term{Paschke dual algebra}.

\begin{definition}
  \label{def:paschke-dual}
  Let $\A,\B$ be separable \cstar-algebras with $\A$ unital and $\B$ stable.
  Let $\phi : \A \to \Mul(\B)$ be a unital absorbing trivial extension.
  The \term{Paschke dual algebra} $\paschkedual{\A}{\B}$ is the relative commutant in the corona algebra $\C(\B)$ of the image of $\A$ under $\pi \circ \phi$, i.e., $\paschkedual{\A}{\B} \defeq \big(\pi \circ \phi (\A)\big)'$.
\end{definition}

The Paschke dual algebra is independent, up to *-isomorphism, of the choice of the absorbing (strongly) unital trivial extension $\phi$.
Indeed, if $\phi, \psi : \A \to \Mul(\B)$ are unital absorbing trivial extensions, then $\phi \sim \phi \oplus \psi \sim \psi$, and so there is a unitary $U \in \Mul(\B)$ such that $\pi(U)(\pi \circ \phi(a))\pi(U^{*}) = \pi \circ \psi(a)$ for all $a \in \A$.
Then this unitary also conjugates the relative commutants of $\pi \circ \phi(\A)$ and $\pi \circ \psi(\A)$.

In \cite{LN-2020-IEOT}, we proved several results about $\paschkedual{\A}{\B}$, generalizing those studied by Paschke for $\paschkedual{\A}{\K}$, which we summarize in the next lemma.

\begin{lemma}[\protect{\cite[Lemma~2.2]{LN-2020-IEOT}}]
  \label{lem:paschke-dual-properly-infinite}
  Let $\A$ be a unital separable nuclear \cstar-algebra, and let $\B$ be a separable stable \cstar-algebra.
  Then we have the following:
  \begin{enumerate}
  \item The Paschke dual algebra $\paschkedual{\A}{\B}$ is properly infinite.  In fact, $1 \oplus 0 \sim 1 \oplus 1$ in $\mat_2 \otimes \paschkedual{\A}{\B}$, i.e., $\paschkedual{\A}{\B}$ contains a unital copy of the Cuntz algebra $\mathcal{O}_2$.
  \item Every element of $\kz(\paschkedual{\A}{\B})$ is represented by a full properly infinite projection in $\paschkedual{\A}{\B}$.
  \end{enumerate}
  \label{lem:Jan1520216AM}
\end{lemma}

In \cite{LN-2020-IEOT}, we also provided the following double commutant theorem for the Paschke dual algebra, which is akin to \cite[Theorem~1]{Ng-2018-NYJM}, and shows that the algebra $\paschkedual{\A}{\B}$ is dual in yet another way.   This generalizes a remark of Valette \cite{Val-1983-PJM}.

\begin{theorem}[\protect{\cite[Theorem~2.10]{LN-2020-IEOT}}]
  \label{thm:standard-position-relative-commutants}
  Let $\A$ be a separable simple unital nuclear \cstar-algebra, and let $\B$ be a separable stable simple \cstar-algebra.
  For any unital trivial absorbing extension $\phi : \A \to \Mul(\B)$,
  \begin{equation*}
    \big(\pi \circ \phi(\A)\big)'' = \pi \circ \phi (\A).
  \end{equation*}
  Equivalently, if we identify $\A$ with its image $\pi \circ \phi (\A)$ (since $\pi \circ \phi$ is injective), and $\paschkedual{\A}{\B} \defeq \big(\pi \circ \phi(\A)\big)'$, then
  \begin{equation*}
    \A' = \paschkedual{\A}{\B} \quad\text{and}\quad (\paschkedual{\A}{\B})' = \A.
  \end{equation*}
\end{theorem}

Towards generalizations of \Cref{thm:bdf-unitary-equivalence-semifinite-factors},
we remind the reader of the following definition due to \cite{DE-2001-KT}.

\begin{definition}
  \label{def:proper-asymptotic-unitary-equivalence}
  Let $\A, \B$ be separable \cstar-algebras, and let
  $\phi, \psi : \A \rightarrow \Mul(\B)$ be *-homomorphisms.

  $\phi$ and $\psi$ are said to be \term{properly asymptotically unitarily equivalent}
  ($\phi \paue \psi$) if there exists a norm continuous path
  $\{ u_t \}_{t \in [0, \infty)}$ of unitaries in $\complex 1 + \B$ such that
  \begin{equation*}
    u_t \phi(a) u_t^* - \psi(a) \in \B, \makebox{   for all }t
  \end{equation*}
  and
  \begin{equation*}
    \norm{ u_t \phi(a) u_t^* - \psi(a) } \rightarrow 0 \makebox{   as  } t \rightarrow
    \infty,
  \end{equation*}
  for all $a \in \A$.
\end{definition}

Note that proper asymptotic unitary equivalence, as above, is a ``local" type of
equivalence since the unitaries are in $\mathbb{C}1 + \B$ and the
two maps $\phi(a) - \psi(a) \in \B$, i.e, $\pi \circ \phi(a) = \pi \circ 
\psi(a)$ for all $a \in \A$.  This is fitting since we will be using this
notion to generalize the BDF essential codimension result, a result about
``local" differences between projections.

Towards stating the next result, we also briefly recall 
the generalized homomorphism picture of $KK = \kkz$, which is also a ``local"
notion. A good reference for this material is \cite[Chapter~4]{JT-1991}
(see also \cite{Bla-1998}).
Let $\A, \B$ be C*-algebras with $\A$ separable and $\B$ 
$\sigma$-unital. A \emph{$KK_h(\A, \B)$-cycle} is a pair $(\phi, \psi)$ of 
*-homomorphisms $\phi, \psi: \A \rightarrow \Mul(\B \otimes \K)$ such that
$\phi(a) - \psi(a) \in \B \otimes \K$ for all $a \in \A$.    
Recall that
two $KK_h(\A, \B)$-cycles $(\phi_0, \psi_0)$ and $(\phi_1, \psi_1)$ 
are \emph{homotopic}  
if there exists a path $\{ (\phi_t, \psi_t) \}_{t \in [0,1]}$ of $KK_h(\A, \B)$-cycles such that for all $a \in \A$,
\begin{enumerate}
\item[i.] the maps $[0,1] \rightarrow \Mul(\B \otimes \K)$ given
by $t \mapsto \phi_t(a)$ and $t \mapsto \psi_t(a)$ are both strictly continuous,
and 
\item[ii.] the map $[0, 1] \rightarrow \B \otimes \K : t \mapsto 
\phi_t(a) - \psi_t(a)$ is norm continuous.
\end{enumerate}
\noindent For any $KK_h(\A, \B)$-cycle $(\phi, \psi)$, we let $[\phi, \psi]$
denote its homotopy class. Then we may define 
$\kkz(\A, \B)$ or $KK(\A, \B)$ by 
$$KK(\A, \B) =_{df} \{ [\phi, \psi] : (\phi, \psi) \makebox{ is a  }
KK_h(\A, \B)-\makebox{cycle} \}.$$
With an appropriate addition operation,
$KK(\A, \B)$ is an abelian group.
For this and 
more information about the bivariant functor $KK$, we refer the
reader to \cite{JT-1991} and \cite{Bla-1998}.  We see already, 
from the definition
as well as the next result, that $KK =\kkz$ is concerned with ``local"
phenomena.
In fact, a very good example of a $\kk_h$-cycle is BDF's definition of essential codimension (see the paragraph after \Cref{thm:bdf-unitary-equivalence-semifinite-factors}).

Recall that a unital C*-algebra $\D$ is said to be \emph{$K_1$-injective} if the standard map
$U(\D)/U_0(\D) \rightarrow K_1(\D)$ is injective. 
As preparation for the next result, we further remark that the corona algebra of a separable stable \cstar-algebra is $\ko$-injective \cite[Proposition~4.9]{GR-2020-GMJ}.

The next result says that $\ko$-injectivity of the Paschke dual $\paschkedual{\A}{\B}$ is sufficient to guarantee interesting uniqueness theorems which generalize \Cref{thm:bdf-unitary-equivalence-semifinite-factors}.  
The proof is already contained in previous works, though often implicitly (\cite{DE-2001-KT,Lee-2011-JFA,LN-2020-IEOT, Lin-2002-JOT}). 
The general strategy goes back to Lin's proof of the first general stable
uniqueness theorem, which was a groundbreaking result (\cite[Theorem~4.3]{Lin-2002-JOT}).
As part of the purpose of this section is to provide an accessible reference,
for the convenience of the reader and to help clean up the literature, we explicitly provide the statement and argument.

\begin{theorem}
  Let $\A$ be a separable nuclear \cstar-algebra, and let $\B$ be a separable stable \cstar-algebra.
  Suppose that either $\A$ is unital and $\paschkedual{\A}{\B}$ is $\ko$-injective, or $\A$ is nonunital and $\paschkedual{(\widetilde{\A})}{\B}$ is $\ko$-injective.
  Let $\phi, \psi : \A \rightarrow \Mul(\B)$ be two absorbing trivial extensions with $\phi(a) - \psi(a) \in \B$ for all $a \in \A$ such that either both $\phi$ and $\psi$ are unital or both $\pi \circ \phi$ and $\pi \circ \psi$ are nonunital.

  Then
  $[\phi, \psi] = 0$ in $\kk(\A, \B)$ if and only if $\phi \paue \psi$.
  \label{thm:MainUniqueness}
\end{theorem}

\begin{proof}

  The ``if'' direction follows directly from Lemma~3.3 of \cite{DE-2001-KT}.

  We now prove the ``only if'' direction.

  Let $\A^+$ denote the unitization of $\A$ if $\A$ is nonunital, and
  $\A \oplus \complex$ if $\A$ is unital.
  Now if $\phi : \A \rightarrow \Mul(\B)$ is a
  nonunital absorbing extension (so $\pi \circ \phi$ is nonunital),
  then $\phi(\A)^{\perp}$ contains a projection which is Murray--von Neumann equivalent
  to $1_{\Mul(\B)}$, and by
  \cite[Section~16, Page~402]{EK-2001-PJM} and \cite{Gab-2016-PJM},
  the map $\phi^+ :
  \A^+ \rightarrow \Mul(\B)$ given by
  $\phi^+ |_{\A} = \phi$ and
  $\phi^+(1) = 1$ is a unital absorbing trivial extension (i.e., $\pi \circ \phi^+$ is unitally absorbing).
  The same holds for $\psi$ and $\psi^+$.
  Moreover, $(\phi^+,\psi^+)$ is a generalized homomorphism.
  Additionally, $[\phi^+,\psi^+] = 0$ because a homotopy of generalized homomorphisms $(\phi_s,\psi_s)$ between $(\phi,\psi)$ and $(0,0)$ lifts to a homotopy $(\phi^+_s,\psi^+_s)$, and hence $[\phi^+,\psi^+] = [0^+,0^+] = 0$.
  Thus, we may assume that
  $\A$ is unital and $\phi$ and $\psi$ are unital *-monomorphisms.

  As before, we may identify the Paschke dual algebra as $\paschkedual{\A}{\B} = (\pi \circ \phi(\A))' \subseteq \C(\B)$.

  By \cite[Lemma~3.3]{LN-2020-IEOT} there exists a norm continuous path  $\{ u_t \}_{t \in [0, \infty)}$ of unitaries in $\Mul(\B)$ such that
  \begin{equation*}
    u_t \phi(a) u_t^* - \psi(a) \in \B
  \end{equation*}
  for all $t$ and for all $a \in \A$, and
  \begin{equation*}
    \norm{u_t \phi(a) u_t^* - \psi(a)} \rightarrow 0
  \end{equation*}
  as $t \rightarrow \infty$, for all $a \in \A$.

  It is trivial to see that this implies that
  \begin{equation*}
    [\phi, u_0\phi u_0^*] = [ \phi, \psi ] = 0,
  \end{equation*}
  and that $\pi(u_t) \in (\pi \circ \phi(\A))' = \paschkedual{\A}{\B}$ for all $t$.

  It is well-known that we have
  a group isomorphism $\kk(\A, \B) \rightarrow \kkhig(\A, \B) : [\phi, \psi]
  \rightarrow [\phi, \psi, 1]$.
  (Here, $\kkhig$ is the version of $\kk$-theory presented in \cite[Section~2]{Hig-1987-PJM}.)
  Hence, $[\phi, u_0 \phi u_0^*, 1] = 0$ in $\kkhig(\A, \B)$.
  Hence, by \cite[Lemma~2.3]{Hig-1987-PJM}, $[\phi, \phi, u_0^*] = 0$
  in $\kkhig(\A, \B)$.

  By Thomsen's Paschke duality theorem \cite[Theorem~3.2]{Tho-2001-PAMS},
  there is a group isomorphism $\ko(\paschkedual{\A}{\B}) \rightarrow \kkhig(\A, \B)$
  which sends $[\pi(u_0)]$ to $[\phi, \phi, u_0^*]$.
  % In general, the map sents $[w] \mapsto [\phi, \phi, w_0]$ where
  % $w_0 \in \pi^{-1}(\paschkedual{\A}{\B})$ and $\pi(w_0) = w$.
  Hence, $[\pi(u_0)] = 0$ in $\ko(\paschkedual{\A}{\B})$.
  Since $\paschkedual{\A}{\B}$ is $\ko$-injective, $\pi(u_0) \sim_h 1$ in $\paschkedual{\A}{\B} =
  (\pi \circ \phi(\A))'$.
  Hence, there exists a unitary $v \in \complex 1 + \B$ such that
  $v^* u_0 \sim_h 1$ in $\pi^{-1}(\paschkedual{\A}{\B})$.
  Indeed, we can construct such a $v$ as follows:
  Since $\pi(u_0) \sim_h 1$ in $\paschkedual{\A}{\B}$, we may write $\pi(u_0)^{*} = e^{i a_1} \cdots e^{i a_n}$
  and hence $e^{-i a_n} \cdots e^{-i a_1} \pi(u_0)^{*} = 1$
  for some self-adjoint $a_1,\ldots, a_n \in \paschkedual{\A}{\B}$.
  Lift $a_1,\ldots,a_n$ to self-adjoint $A_1,\ldots,A_n \in \pi^{-1}(\paschkedual{\A}{\B})$.
  Then set $v = e^{-i A_n} \cdots e^{-i A_1} u_0^{*} \in 1 + \B$.

  Hence, modifying an initial segment of $\{ v^*u_t \}_{t \in [0, \infty)}$
  if necessary, we may assume that $\{ v^*u_t \}_{t \in [0, \infty)}$ is
  a norm continuous path of unitaries in $\pi^{-1}(\paschkedual{\A}{\B})$ such that
  $v^*u_0 = 1$.

  Now for all $t \in [0,\infty)$, let $\alpha_t \in \Aut(\phi(\A) + \B)$ be given
  by $\alpha_t(x) \defeq v^* u_t x u_t^* v$ for all $x \in \phi(\A) + \B$.
  Thus, $\{ \alpha_t \}_{t \in [0, \infty)}$ is a uniformly continuous path of
  automorphisms of $\phi(\A) + \B$ such that $\alpha_0 = \id$.
  Hence, by \cite[Proposition~2.15]{DE-2001-KT} (see also
  \cite[Theorems~3.2 and 3.4]{Lin-2002-JOT}),
  there exists a continuous path $\{ v_t \}_{t\in [0, \infty)}$ of unitaries
  in $\phi(\A) + \B$ such that
  $v_0 = 1$ and
  $\norm{ v_t x v_t^* - v^* u_t x u_t^* v } \rightarrow 0$ as $t \rightarrow \infty$
  for all $x \in \phi(\A) + \B$.
  Thus, $\norm{ v v_t x v_t^* v^* - u_t x u_t^* } \rightarrow 0$ as
  $t \rightarrow \infty$
  for all $x \in \phi(\A) + \B$.

  We now proceed as in the last part of the proof of \cite[Proposition 3.6. Step 1]{DE-2001-KT} (see
  also the proof of \cite[Theorem~3.4]{Lin-2002-JOT}).
  For all $t \in [0, \infty)$, let $a_t \in \A$ and $b_t \in \B$ such that
  $v v_t = \phi(a_t) + b_t$.
  Since $\pi \circ \phi$ is injective, we have that for all $t$,
  $a_t$ is a unitary in $\A$, and hence, $\phi(a_t)$ is a unitary in $\phi(\A) + \B$.
  Note also that since $\pi \circ \phi = \pi \circ \psi$ and both maps are
  injective, $\norm{ a_t a a_t^* - a } \rightarrow 0$ as $t \rightarrow \infty$ for
  all $a \in \A$.
  For all $t$, let $w_t \defeq v v_t \phi(a_t)^* \in 1 + \B$.
  Then $\{ w_t \}_{t \in [0,1)}$ is a norm continuous path of unitaries in
  $\complex 1 + \B$, and for all $a \in \A$,
  \begin{align*}
    \norm{w_t \phi(a) w_t^* - \psi(a)}
    &\leq \norm{ w_t \phi(a) w_t^* - v v_t \phi(a) v_t^* v^* } \\
    &\qquad+ \norm{ v v_t \phi(a) v_t^* v^* - u_t \phi(a) u_t^* } \\
    &\qquad+ \norm{ u_t \phi(a) u_t^* - \psi(a) } \\
    &= \norm{ v v_t\phi( a_t^* a a_t - a) v_t^* v^* } \\
    &\qquad+ \norm{ v v_t \phi(a) v_t^* v^* - u_t \phi(a) u_t^* } \\
    &\qquad+ \norm{ u_t \phi(a) u_t^* - \psi(a) } \\
    &\rightarrow 0. \qedhere
  \end{align*}

\end{proof}

\begin{lemma}[\protect{\cite[Lemma~2.4]{LN-2020-IEOT}}]
  Let $\A$ be a unital separable nuclear \cstar-algebra and $\B$ be a separable
  stable \cstar-algebra such that either $\B \cong \K$ or $\B$ is
  simple purely infinite.

  Then $\paschkedual{\A}{\B}$ is $\ko$-injective.
  \label{lem:PIK1Injective}
\end{lemma}

Note that for a separable simple stable C*-algebra $\B$, $\C(\B)$ is simple if and
only if $\C(\B)$ is simple purely infinite if and only if $\B \cong \K$ or $\B$ is
simple purely infinite \cite[Theorem~3.8]{Lin-1991-PAMS}.   Hence,
the above lemma resolves the problem of $K_1$-injectivity for the Paschke dual
$\A^d_{\B}$ in the case where the corona algebra $\C(\B)$ is simple.

\begin{theorem}
  Let $\A$, $\B$ be separable \cstar-algebras such that $\A$ is nuclear and
  either $\B \cong \K$ or $\B$ is stable simple purely infinite.
  Let $\phi, \psi : \A \rightarrow \Mul(\B)$ be essential trivial extensions
  with $\phi(a) - \psi(a) \in \B$ for all $a \in \A$  such that
  either both $\phi$ and $\psi$ are unital or both $\pi \circ \phi$ and
  $\pi \circ \psi$ are nonunital.

  Then $[\phi, \psi] = 0$ in $\kk(\A, \B)$ if and only if $\phi \paue \psi$.
\end{theorem}

\begin{proof}
  This follows immediately from \Cref{thm:MainUniqueness} and
  \Cref{lem:PIK1Injective}.

  Of course, we here are using that we are in the ``nicest'' setting for
  extension theory:  By our hypotheses on $\A$ and $\B$, every essential
  extension of $\B$ by $\A$ is absorbing in the appropriate sense (see, for example, \cite[Proposition~2.5]{GN-2019-AM} and \cite[Theorem~17]{EK-2001-PJM}).
\end{proof}

Thus, from the point of view of simple stable canonical ideals with appropriate
regularity properties,
what remains is the case where the canonical ideal is stably finite.  In
\cite{LN-2020-IEOT}, we had a partial result with restrictive conditions
(see the present paper \Cref{thm:Jan152021Restrictive}).  Part of
the goal of the present paper is to remove many of these restrictive conditions
(see the present paper \Cref{thm:main-theorem}).

The next result partially answers a question that Professor Huaxin Lin
asked the second author.  The argument actually comes from Lin's
paper \cite[Proposition~2.6]{Lin-2005-CJM}.

\begin{proposition}
  \label{prop:paschke-dual-purely-infinite}
  Let $\A$ be a unital separable simple nuclear \cstar-algebra, and let
  $\B$ be a separable stable \cstar-algebra such that
  either
  \begin{equation*}
    \B \cong \K \makebox{  or  } \B \makebox{  is simple purely infinite.}
  \end{equation*}

  Let $\sigma : \A \rightarrow \Mul(\B)$ be a unital trivial essential
  extension.

  Then $(\pi \circ \sigma (\A))'$  is simple purely infinite.
  As a consequence, $\paschkedual{\A}{\B}$ is simple purely infinite.
\end{proposition}

\begin{proof}
  Note that $\sigma$ is absorbing (e.g., see
  \cite[Proposition~2.5]{GN-2019-AM}).

  By \Cref{lem:paschke-dual-properly-infinite}, $\paschkedual{\A}{\B}$ contains a unital copy of $\mathcal{O}_2$, so it cannot be isomorphic to $\complex$.
  Let $c \in (\pi \circ \sigma(\A))'_+ - \{ 0 \}$ be arbitrary.
  We want to find $s \in (\pi \circ \sigma(\A))'$ so that
  $s c s^* = 1$.

  We may assume that $\norm{c} = 1$.

  Let
  \begin{equation*}
    X \defeq \spec(c).
  \end{equation*}

  Let
  \begin{equation*}
    \phi, \psi : C(X) \otimes \A \rightarrow \Mul(\B)/\B
  \end{equation*} be *-homomorphisms
  that are given as follows (using the universal property of the maximal tensor product):

  \begin{equation*}
    \phi : f \otimes a \mapsto f(c) (\pi \circ \sigma(a))
  \end{equation*}
  and
  \begin{equation*}
    \psi : f \otimes a \mapsto f(1) (\pi \circ \sigma(a)).
  \end{equation*}

  Both $\phi$ and $\psi$ are unital extensions.  Since $\A$ is simple unital, it is a short exercise to see that $\phi$ is essential.  Moreover,
  one can check that $\psi$ is a strongly unital trivial extension.

  Hence, since either $\B \cong \K$ or $\B$ is simple purely infinite,
  we must have that $\phi$ is absorbing (from \cite[Proposition~2.5]{GN-2019-AM}), so
  \begin{equation*}
    \phi \sim \phi \oplus \psi
  \end{equation*}
  where $\oplus$ is the BDF sum, and $\sim$ is unitary equivalence with unitary coming from $\Mul(\B)$.

  So let $W \in \mat_2 \otimes \Mul(\B)$ be such that
  \begin{equation*}
    W^* W = 1 \oplus 1,
  \end{equation*}
  \begin{equation*}
    W W^* = 1 \oplus 0,
  \end{equation*}
  and letting $w \defeq \pi(W)$
  \begin{equation*}
    w^*(\phi(\cdot) \oplus 0)w = \psi(\cdot) \oplus \phi(\cdot).
  \end{equation*}

  Then
  \begin{align*}
    w^*(c \oplus 0)w
    &= w^*(\phi(\id_X \otimes 1) \oplus 0)w\\
    &= \psi(\id_X \otimes 1) \oplus \phi(\id_X \otimes 1)\\
    &= 1 \oplus c. 
  \end{align*}

  So
  \begin{equation*}
    \left[\begin{array}{cc} 1 & 0 \\ 0 & 0 \end{array}\right]
    w^* (c \oplus 0) w \left[\begin{array}{cc} 1 & 0 \\ 0 & 0 \end{array} \right]
    = 1 \oplus 0.
  \end{equation*}

  Identifying $\Mul(\B)$ with $e_{1,1} \otimes \Mul(\B)$, we may view
  $w \left[ \begin{array}{cc} 1 & 0 \\ 0 & 0 \end{array} \right]$ as an
  element of $\C(\B)$.  Hence, from the above computation, to finish the
  proof, it suffices to prove that
  \begin{equation*}
    w \left[\begin{array}{cc} 1 & 0 \\ 0 & 0 \end{array} \right]
    \in (\pi \circ \sigma(\A))'.
  \end{equation*}

  But for all $a \in \A$,
  \begin{align*}
    w^* ( \pi \circ \sigma(a) \oplus 0) w
    &=  w^* (\phi(1 \otimes a) \oplus 0)w \\
    &= \psi(1 \otimes a) \oplus \phi(1 \otimes a)\\
    &= \pi \circ \sigma(a) \oplus \pi \circ \sigma(a).
  \end{align*}

  So
  \begin{equation*}
    (\pi \circ \sigma(a) \oplus 0) w =
    w ( \pi \circ \sigma(a) \oplus \pi \circ \sigma(a))
  \end{equation*}
  for all $a \in \A$.

  Since $ww^* = 1 \oplus 0$,
  \begin{equation*}
    (\pi \circ \sigma(a) \oplus \pi \circ \sigma(a))w =
    w (\pi \circ \sigma(a) \oplus \pi \circ \sigma(a))
  \end{equation*}
  for all $a \in \A$.

  Hence,
  \begin{equation*}
    w \in \mat_2 \otimes (\pi \circ \sigma(\A))'.
  \end{equation*}

  Hence,
  \begin{equation*}
    w \left[ \begin{array}{cc} 1 & 0 \\ 0 & 0 \end{array} \right]
    \in (\pi \circ \sigma(\A))'
  \end{equation*}
  as required.
\end{proof}

Our next result establishes a generic sufficient criterion for establishing the $\ko$-injectivity of the Paschke dual algebra.
It is the key technique used in the next section.

Let $\C \subseteq \D$ be an inclusion of \cstar-algebras. Recall that
$\C$ is said to be \term{strongly full} in $\D$ if
every nonzero element of $\C$ is full in $\D$, i.e., for all $y \in \C - \{ 0 \}$, $\Ideal_{\D}(y) = \D$.
Recall also that an element $x \in \D$ is said to be \term{strongly full} in $\D$ if $C^*(x)$ is strongly full in $\D$.
Finally, an extension $\phi : \A \to \C(\B)$ is said to be \term{full} if $\phi$ is essential and $\phi(\A)$ is strongly full in $\C(\B)$.

We remark on a key fact we use in the proof: If $\A$ is a unital separable nuclear \cstar-algebra and $\B$ is a separable stable with the corona factorization property (CFP), then every unital full extension of $\B$ by $\A$ is absorbing \cite{KN-2006-HJM}.
(In fact, this is a characterization of the corona factorization property.)
Many separable simple stable \cstar-algebras have the corona factorization property, including all such \cstar-algebras with strict comparison of positive elements, and also all such \cstar-algebras which are purely infinite (see \cite[Theorems~5.11 and 5.13, Definition~2.12]{OPR-2012-IMRN} and \cite{OPR-2011-TAMS}).
In particular, for any separable simple and purely infinite \cstar-algebra $\B$, the Cuntz semigroup of $\B$ is just $\{0, \infty\}$, and so it is particularly easy to check that it satisfies \cite[Definition~2.12]{OPR-2012-IMRN}.
For the case where $\B$ is a separable simple \cstar-algebra with strict comparison, note that every $\tau \in T(\B)$ extends to a dimension function on the Cuntz semigroup of $\B$ (see the discussion after \cite[Theorem~4.5]{Ror-2004-IJM}, see also the conclusion of \cite[Corollary~4.7]{Ror-2004-IJM} which gives a characterization of strict comparison in terms of dimension functions on Cuntz semigroups).

\begin{theorem}
  \label{thm:k1-injectivity-generic}
  Suppose that $\A,\B$ are separable \cstar-algebras with $\A$ unital, simple and nuclear, and $\B$ stable and having the corona factorization property.

  Suppose that $\phi : \A \to \Mul(\B)$ is a unital trivial absorbing extension and realize $\paschkedual{\A}{\B}$ as $(\pi \circ \phi (\A))'$.
  If every unitary in the Paschke dual algebra $\paschkedual{\A}{\B}$ is homotopic in $\paschkedual{\A}{\B}$ to a unitary which is strongly full in $\C(\B)$, then $\paschkedual{\A}{\B}$ is $\ko$-injective.
  Moreover, for each $n \in \nat$, the map
  \begin{equation*}
    U(\mat_n \otimes \paschkedual{\A}{\B})/U(\mat_n \otimes \paschkedual{\A}{\B})_0 \to U(\mat_{2n} \otimes \paschkedual{\A}{\B})/U(\mat_{2n} \otimes \paschkedual{\A}{\B})_0
  \end{equation*}
  given by
  \begin{equation*} [u] \mapsto [u \oplus 1] \end{equation*}
  is injective.
\end{theorem}

\begin{proof}  
  Clearly, $\ko$-injectivity will follow from the more specific statement.
  By \Cref{lem:paschke-dual-properly-infinite}, $\paschkedual{\A}{\B}$ contains a unital copy of $\mathcal{O}_2$, and so $1 \oplus 1 \sim 1$.
  Therefore, for all $n \in \nat$, $\paschkedual{\A}{\B} \cong \mat_n \otimes \paschkedual{\A}{\B}$.
  Hence it suffices to establish injectivity for the map in the $n = 1$ case: \begin{equation*}
    U(\paschkedual{\A}{\B})/U(\paschkedual{\A}{\B})_0 \to U(\mat_2 \otimes \paschkedual{\A}{\B})/U(\mat_2 \otimes \paschkedual{\A}{\B})_0.
  \end{equation*}

  Let $u \in \paschkedual{\A}{\B}$ be a unitary for which
  \begin{equation*}
    u \oplus 1 \sim_h 1 \oplus 1
  \end{equation*} in $\mat_2 \otimes \paschkedual{\A}{\B}$.

  By hypothesis, $u$ is homotopic, in $\paschkedual{\A}{\B}$, to a unitary which is strongly full in $\C(\B)$.
  Therefore, we may assume, without loss of generality, that $u$ is strongly full in $\C(\B)$.
  By \cite[Lemma~2.6]{LN-2020-IEOT}, $C^{*}(\pi \circ \phi(\A),u)$ is strongly full in $\C(\B)$.
  Hence, the inclusion map
  \begin{equation*}
    \i : C^{*}(\pi \circ \phi(\A),u) \hookrightarrow \C(\B)
  \end{equation*}
  is a full extension.  Since $\B$ has the corona factorization property, and since $C^*(\pi \circ \phi(\A), u)$ is nuclear (being a quotient of $C(\mathbb{S}^1) \otimes \A$), $\i$ is a unital absorbing extension.

  Let $\sigma : C^{*}(\pi \circ \phi(\A),u) \to \Mul(\B)$ be a unital trivial absorbing extension.
  Since the restriction $\sigma |_{\pi \circ \phi(\A)}$ is also a unital trivial absorbing extension, conjugating
  $\sigma$ by a unitary if necessary, we may assume that the restriction of $\pi \circ \sigma$ to $\pi \circ \phi(\A)$ is the identity
  map.
  (See the paragraph after \Cref{def:paschke-dual} for why we have unitary equivalence.)

  By \cite[Lemma 2.3]{LN-2020-IEOT},
  \begin{equation*}
    \pi \circ \sigma(u) \sim_h 1
  \end{equation*}
  in $\paschkedual{\A}{\B}$.

  Since $\i$ is absorbing,
  \begin{equation*}
    \i \oplus (\pi \circ \sigma) \sim \i.
  \end{equation*}

  Consequently, there is an isometry $\tilde{v} \in \mat_2 \otimes \Mul(\B)$ such that $\tilde{v}^{*} \tilde{v} = 1 \oplus 1$, $\tilde{v} \tilde{v}^{*} = 1 \oplus 0$, and for $v \defeq \pi(\tilde{v})$,
  \begin{equation*}
    v (\i \oplus (\pi \circ \sigma)) v^{*} = \i \oplus 0.
  \end{equation*}
  Therefore,
  \begin{equation*}
    v(a \oplus a)v^{*} = a \oplus 0
  \end{equation*} for all $a \in \pi \circ \phi(\A)$, and also
  \begin{equation*}
    v(u \oplus (\pi \circ \sigma(u)))v^{*} = u \oplus 0.
  \end{equation*}
  Thus, since $vv^{*} = 1 \oplus 0$,
  \begin{equation*}
    v(a \oplus a) = (a \oplus 0)v = (a \oplus a) v.
  \end{equation*} for all $a \in \pi \circ \phi(\A)$, and therefore,
  \begin{equation*}
    v \in \mat_2 \otimes \paschkedual{\A}{\B}.
  \end{equation*}

  Moreover, we also have
  \begin{equation*}
    u \oplus (\pi \circ \sigma (u)) \sim_h u \oplus 1 \sim_h 1 \oplus 1
  \end{equation*}
  within $\mat_2 \otimes \paschkedual{\A}{\B}$.
  Conjugating this continuous path of unitaries by $v$, we obtain
  \begin{equation*}
    u \sim_h 1
  \end{equation*} within $\paschkedual{\A}{\B}$.
\end{proof}

\section{The Paschke dual for simple nuclear \cstar-algebras}
\label{sec:paschke-dual-simple}

We first fix some notation, which will be used for the rest of this paper.
For a $\sigma$-unital simple \cstar-algebra $\B$ and for a nonzero element $e \in \Ped(\B)_+$
(the Pedersen ideal of $\B$), we let
$T_e(\B)$ denote the set of all densely defined, norm-lower semicontinuous traces
$\tau$ on $\B_+$ such that $\tau(e) = 1$.  Recall that $T_e(\B)$, with the topology
of pointwise convergence on $\Ped(\B)$, is a compact convex set.  In fact, it is
a Choquet simplex (see \cite{ERS-2011-AJM, TT-2015-CMB}).  All the results and arguments in this paper are independent of
the choice of normalizing element $e \in \Ped(\B)_+ - \{ 0 \}$.  Hence, we will
usually drop the $e$ and simply write $T(\B)$.
Recall also, that every $\tau \in T(\B)$ extends uniquely to a
strictly lower semicontinuous trace
$\Mul(\B)_+ \rightarrow [0, \infty]$, which we will also denote by ``$\tau$''.

Recall next that for all
$\tau \in T(\B)$ and for any $a \in \Mul(\B)_+$,
\begin{equation*}
  d_{\tau}(a) \defeq \lim_{n\rightarrow \infty} \tau(a^{1/n}) \in [0, \infty].
\end{equation*}

For $\delta > 0$,
let $f_{\delta} : [0, \infty) \rightarrow [0,1]$ be the unique continuous function for which
\begin{equation*}
  f_{\delta}(t)
  =
  \begin{cases}
    1 & t \in [\delta, \infty) \\
    0 & t = 0\\
    \makebox{linear on  } & [0, \delta].
  \end{cases}
\end{equation*}

In what follows, for elements $a,b$ in a \cstar-algebra, $a \approx_{\epsilon} b$ means $\norm{ a-b } < \epsilon$.

\subsection{Homotopy trivial unitaries in a Paschke dual algebra}

The  lemmas in this section give a technique for constructing unitaries, in a  
Paschke dual algebra, which are homotopy equivalent to $1$.\\

\begin{df}  Let $\{ \alpha_n \}$ be a sequence of complex numbers.
$\{ \alpha_n \}$ is said to have the \emph{close neighbors property} if
$$\lim_{n \rightarrow \infty } |\alpha_{n+1} - \alpha_n | = 0,$$ 
i.e., 
for every $\epsilon > 0$, there exists an $N \geq 1$ such that for all $n \geq N$,
$$|\alpha_{n+1} - \alpha_n | < \epsilon.$$
\label{df:CloseNeighborsProperty}
\end{df}

\begin{lemma}
 
 Let $\B$ be a nonunital $\sigma$-unital simple \cstar-algebra.

  Suppose that
  $\{ e_n \}_{n=1}^{\infty}$ is an approximate unit for
  $\B$ such that
  \begin{equation*}
    e_{n+1} e_n = e_n \makebox{  for all  } n\geq 1.
  \end{equation*}

  Suppose that $\{ \alpha_n \}_{n=1}^{\infty}$ is a sequence in
  the unit circle (of the complex plane) $S^1$ with the close neighbors property. 

  Then $\pi\left( \sum_{n=1}^{\infty} \alpha_n (e_n - e_{n-1}) \right)$
  is a unitary in $\C(\B)$.
  (We use the convention $e_0 \defeq 0$.)

\label{lem:CloseNeighborsUnitary}
\label{lem:AUUnitary}
\end{lemma}

\begin{proof}
For all $n \geq 1$, let $r_n =_{df} e_n - e_{n-1}$.
Let $X \in \Mul(\B)$ be given by 
$$X =_{df} \sum_{n=1}^{\infty} \alpha_n r_n,$$
where the sum converges strictly in $\Mul(\B)$.   

Note that since $e_{n+1} e_n = e_n$ for all $n \geq 1$,
$r_m \perp r_n$ for all $|m - n | \geq 2$.

Also note that by the hypotheses on $\{ \alpha_n \}$,
$\overline{\alpha_n} \alpha_{n+1} \rightarrow 1$ as $n \rightarrow
\infty$.

Hence,
\begin{eqnarray*}
X^* X & = & \left(\sum_{n=1}^{\infty} \overline{\alpha_n} r_n \right)
\left(\sum_{k=1}^{\infty} \alpha_k r_k \right) \\
& = & \sum_{n=1}^{\infty} \overline{\alpha_n} r_n (\alpha_{n-1} r_{n-1}
+ \alpha_n r_n + \alpha_{n+1} r_{n+1}) \makebox{   (Note  } e_0 =_{df}
e_{-1} =_{df} 0 \makebox{)} \\
& = & \sum_{n=1}^{\infty} r_n (\overline{\alpha_n}\alpha_{n-1} r_{n-1}
+ r_n + \overline{\alpha_n}\alpha_{n+1} r_{n+1}) \\
& = & \sum_{n=1}^{\infty} r_n (r_{n-1} + r_n + r_{n+1}) \makebox{  (mod }
\B \makebox{)}\\
& = & \left(\sum_{n=1}^{\infty} r_n \right)\left(\sum_{k=1}^{\infty} r_k
\right)\\
& = & (1_{\Mul(\B)})(1_{\Mul(\B)})\\
& = &  1_{\Mul(\B)}. 
\end{eqnarray*}  

By a similar argument, $X X^* = 1_{\Mul(\B)}$ (mod $\B$).

So $\pi(X)$ is a unitary in $\C(\B)$.   
\end{proof}

The next computation should be well-known, but we provide it for the convenience of the
reader.

Let $\J$ be an ideal of a C*-algebra $\E$, let $S \subseteq \E$ and
let $\{ e_{\gamma} \}$ be an approximate unit for $\J$. Recall
that  $\{ e_{\gamma} \}$ \emph{quasicentralizes} $S$ if
for all $x \in S$, $\| e_{\gamma} x - x e_{\gamma} \| \rightarrow 0$.

\begin{lemma}
  \label{lem:QuasicentralAU}
  Let $\B$ be a nonunital $\sigma$-unital simple \cstar-algebra, and let
  $S \defeq \{ x_k \mid k \geq 1 \} \subseteq \Mul(\B)$
  be a countable set.

  Then there exists an approximate unit $\{ e_n \}$ for $\B$ such that
  \begin{equation*}
    e_{n+1} e_n = e_n
  \end{equation*}
  for all $n$,
  and
  $\{ e_n \}$ quasicentralizes $S$ in the following strong sense:
  For all $n$, for all $1 \leq k \leq n$,
  \begin{equation*}
    \norm{ e_n x_k - x_k e_n } < \frac{1}{2^n}.
  \end{equation*}

  Moreover, for any $\{ e_n \}$ as above, we have that
  \begin{equation*}
    \pi\left(\sum_{n=1}^{\infty} \alpha_n (e_n - e_{n-1})\right) \in \pi(S)'
  \end{equation*}
  for every bounded sequence $\{ \alpha_n \}$ of complex numbers.
\end{lemma}

\begin{proof}

  Let $S \defeq \{ x_k : k \geq 1 \}$.

  Let $\{ e_n \}$ be an approximate unit for $\B$ such that
  \begin{equation*}
    e_{n+1} e_n = e_n
  \end{equation*}
  for all $n$, and $\{ e_n \}$ quasicentralizes $\{ x_k \}$ (see \cite[Theorem~2.2]{Ped-1990});
  that is, for all $k$,
  \begin{equation*}
    \norm{ e_n x_k - x_k e_n } \rightarrow 0
  \end{equation*}
  as $n \rightarrow \infty$.

  Passing to a subsequence of $\{ e_n \}$ if necessary, we may assume that for all $n$, for all $1 \leq k \leq n$,
  \begin{equation*}
    \norm{ e_n x_k - x_k e_n } < \frac{1}{2^n}.
  \end{equation*}

  For all $n$, let
  \begin{equation*}
    r_n \defeq e_n - e_{n-1}.
  \end{equation*}

  Then, for each $k \in \nat$, let $y \defeq \sum_{n=1}^{\infty} \alpha_n r_n$,
where the sum converges strictly in $\Mul(\B)$, and notice that
  \begin{equation*}
    x_k y - y x_k = \sum_{n=1}^{\infty} \alpha_n \big( (x_k e_n - e_n x_k) + (e_{n-1} x_k - x_k e_{n-1} ) \big).
  \end{equation*}
  But for $n > k$,
  \begin{equation*}
    \norm{(x_k e_n - e_n x_k) + (e_{n-1} x_k - x_k e_{n-1} )} \le \frac{3}{2^n},
  \end{equation*}
  and the sequence $\alpha_n$ is bounded.
  Hence the sum representing $x_k y - y x_k$ converges in norm, so $x_k y - y x_k \in \B$.
  I.e., $x_k \left(\sum_{n=1}^{\infty} \alpha_n r_n \right) - \left(\sum_{n=1}^{\infty} \alpha_n r_n \right) x_k \in \B$,
  and therefore $\pi(y) =_{df} \pi\left(\sum_{n=1}^{\infty} \alpha_n (e_n - e_{n-1})\right)$ commutes with each $\pi(x_k)$.
\end{proof}

Morally, the following definition corresponds to a sequence in the unit circle which, when adjacent terms of the sequence are interpolated by the arc between them, corresponds to a path on the unit circle with winding number zero, and therefore corresponds to a unitary which is trivial in $\ko(C(\mathbb{T}))$.

\begin{df}
Let $\{ \alpha_n \}$ be a sequence in the unit circle $S^1$.
Then $\{ \alpha_n \}$ is said to be  \emph{unit oscillating}        
if there exists two subsequences $\{ M_k \}$, $\{ N_k \}$, of the positive
integers such that for all $k$,
\begin{enumerate}
\item[(a)] $M_1 = 1$,  
\item[(b)]  $M_k + 2 < N_k < N_k + 3 < M_{k+1}$,          
\item[(c)]  $\alpha_{M_k - 1} = \alpha_{M_{k}} = 1$ for $k > 1$, 
\item[(d)]  for all $M_k < m < n \leq N_k$,
there exist $0 < r < s \leq 2\pi$ such that
$\alpha_m = e^{ir}$ and $\alpha_n = e^{is}$, and
\item[(e)]  for all $N_k \leq m < n < M_{k+1}-1$,
there exist $0 < r < s \leq 2 \pi$ such that
$\alpha_m = e^{is}$  and   $\alpha_n = e^{ir}$. 
\end{enumerate}
\label{df:UnitOscillatingSequence}
\end{df}

Notice that the ``path" of a unit oscillating sequence never ``strictly" crosses
$1$. It can reach $1$ only at a step of the form $M_k - 1$ or a step of the form
$N_k$.   
Every time that it finally
 ``lands on" $1$ at some step $M_k - 1$, it stays on $1$ for one more step (at step
$M_{k}$), and then it moves in the opposite direction in the
next step (at step $M_k + 1$).  Every time that it finally ``lands on" $1$ at some
step $N_k$, it moves in the opposite direction in the next step (at step $N_k + 1$).  

\begin{lemma}
\label{lem:HomotopyTrivialInCommutant}
Let $\B$ be a nonunital $\sigma$-unital simple C*-algebra, let $S \subset \Mul(\B)$
be a countable set.
Let $\{ e_n \}$ be an approximate unit for $\B$, as in Lemma \ref{lem:QuasicentralAU}
for the given set $S$. 

Let $\{ \alpha_n \}$ be a sequence in $S^1$ which is unit oscillating and having
the close neighbors property, and let
$u =_{df} \pi(\sum \alpha_n (e_n - e_{n-1})) \in \pi(S)'$ be the unitary given
by Lemma \ref{lem:QuasicentralAU}.  

Then $u$ is homotopy trivial in $\pi(S)'$. I.e., there exists a norm-continuous
path $\{ u(t) \}_{t \in [0,2\pi]}$, of unitaries in $\pi(S)'$, such that 
$u(0) = u$ and $u(2\pi) = 1_{\C(\B)}$.
\end{lemma}
           
\begin{proof}
For all $t \in [0,2\pi]$, we define a map
$\theta(t) :[0, 2\pi] \rightarrow [0, 2 \pi]$ 
by $$\theta(t)(r) = \max\{ r - t, 0 \}$$ for all $r \in [0, 2\pi]$.
So $\theta(t)$ is subtraction by $t$, but with a floor of $0$ 
(disallowing negative values).

For all $n$ and all $t \in [0, 2\pi]$, we let 
$\alpha_n(t)$ be defined as follows:

\begin{enumerate}
\item[i.] For all $k$, $\alpha_{M_k - 1}(t) = \alpha_{M_k}(t) =_{df} 1$. 
\item[ii.] For all $k$ and all $M_k < m < M_{k+1} - 1$,
if $\alpha_m = e^{ir}$ and $r \in (0, 2\pi]$, then
$\alpha_m(t)=_{df} exp(i \theta(t)(r))$.  
\end{enumerate}
  
By the definitions of $t \mapsto \theta(t)$ and $t \mapsto \alpha_n(t)$ for all $n$,
we can check that for all $n$, for all $s, t \in [0, 2\pi]$, there is an $r \in [0, 
2\pi]$ such that
$$|\alpha_n(s) - \alpha_n(t)| \leq |e^{i(r- s)} - e^{i(r - t)}|
= |e^{-is} - e^{-it}| \leq 2 |s - t|.$$
So $\{ \alpha_n \}$ is a uniformly equicontinuous family of functions
on $[0, 2\pi]$.  So the map $[0, 2\pi] \rightarrow \Mul(\B) :
t \mapsto \sum \alpha_n(t) (e_n - e_{n-1})$ is norm continuous.                     

Also, since $\{ \alpha_n \}$ is unit oscillating (and by the definitions of
$t \mapsto \theta(t)$ and $t \mapsto \alpha_n(t)$ for all $n$), we see that
for all $t \in [0,2\pi]$, for all $n$,
$$|\alpha_{n+1}(t) - \alpha_n(t)| \leq |\alpha_{n+1} - \alpha_n|.$$
(Recall that a unit oscillating sequence never ``strictly" crosses $1$. Every time
it lands on $1$ (at some $M_k - 1$ for $k \geq 2$ or some $N_k$), 
it either stays on $1$ for one more step and then
    moves in the opposite direction in the next step, or it immediately moves in
the opposite direction in the next step.)                  
Hence, since $\{ \alpha_n \}$ has the close neighbors property, 
for all $t \in [0, 2 \pi]$, $\{ \alpha_n(t) \}$ is a sequence in $S^1$ with
the close neighbors property.  Hence, by Lemma \ref{lem:CloseNeighborsUnitary},
for all $t \in [0, 2 \pi]$, $\pi(\sum \alpha_n(t) (e_n - e_{n-1}))$ is a 
unitary in $\C(\B)$.  Clearly, $u(2 \pi) =1$.   

Finally, by our choice of $\{ e_n \}$ and by Lemma \ref{lem:QuasicentralAU},
for all $t \in [0, 2\pi]$, $\pi(\sum \alpha_n(t)(e_m - e_{n-1})) \in \pi(S)'$.

For all $n$ and all $t \in [0, 2\pi]$, define 
$$u(t) =_{df} \pi\left(\sum \alpha_n(t)(e_m - e_{n-1})\right).$$
\end{proof}

\subsection{Asymptotic properties for Laurent polynomials}

\label{subsec:AsympLaurent}

In this section, we define a series of properties for Laurent polynomials, 
which will facilitate our computations.  We will then use these properties
to prove \Cref{thm:Fullness} and \Cref{thm:Mar2520225PM}.

We define a \term{Laurent polynomial} on the punctured closed
disk  $\overline{B(0,1)} - \{ 0 \}$ to be  a
continuous function $h : \overline{B(0,1)} - \{ 0 \} \rightarrow \complex$ which has the form
\begin{equation*}
  h(\lambda) = \sum_{n= 0}^N \beta_n \lambda^n  + \sum_{m=1}^M \gamma_m \overline{\lambda}^m
\end{equation*}
for all $\lambda \in \overline{B(0,1)} - \{ 0 \}$.
Here, $M, N \geq 1$ are integers and $\beta_n,\gamma_m \in \complex$.

In what follows, we will use that the vector space of Laurent polynomials, when
restricted to the circle $S^1$, is uniformly dense in
$C(S^1)$.

Next, if $\C$ is a unital \cstar-algebra and $h$ is a Laurent polynomial as
above, then for all contractive $x \in \C$, we define
\begin{equation*}
  h(x) \defeq \sum_{n=0}^N \beta_n x^n + \sum_{m=1}^M \gamma_m (x^*)^m.
\end{equation*}
This is well-defined by the uniqueness of Laurent series expansion.  Note that
when $x$ is a unitary, this is consistent with the continuous functional
calculus.

For a real-valued function $f$, we let $\osupp(f) \defeq f^{-1}(\real\setminus\{0\})$ denote the \term{open support} of $f$.
Of course, $\overline{\osupp(f)} = \supp(f)$.\\  

\textbf{
For the rest of this section, we make the following assumptions:
Let $\B$ be a nonunital simple $\sigma$-unital C*-algebra, and let $\Gamma$ be \emph{any}
 collection
of strictly converging series of the form $\sum_{j=1}^{\infty} b_j$ in $\Mul(\B)$, 
where 
$b_j \in \B$ for all $j$.  
(We stress that $\Gamma$ is NOT the set of all such series, but just 
\emph{any} subcollection.) 
Let $\{ e_n \}$ be an approximate unit for $\B$ such that $e_{n+1} e_n = e_n$ for
all $n$.}\\

For the proof of Theorem \ref{thm:Fullness}, we will only need Properties C and X
from the following
definitions.   However, in practice, Properties A and B are used to prove Property
X.  Also, the proof of Property C is similar to those of Properties A and B, and
moreover, Properties A and B are themselves quite interesting.\\ 

\begin{definition}
  \label{def:Properties-A-B-C} \hfill
  \begin{enumerate}
  \item We say that $\Gamma$   has \emph{Property A} with respect to $\{ e_n \}$ 
    if
    for every Laurent polynomial $h$, 
    every contractive $X \in \Mul(\B)$,  every $\epsilon > 0$ and every
    $M \geq 1$, there exists an $L_1 \geq 1$ where for all
    $\sum b_j$ and $\sum b'_j$  in $\Gamma$ for which
    $$b_j = b'_j \makebox{  for all  } j \geq M,$$
    we have that 
    $$h\left(\sum b_j X \right)(1 - e_{L_1}) \approx_{\epsilon} 
    h\left(\sum b'_j X \right)(1 - e_{L_1}).$$\\  

  \item We say that $\Gamma$ 
    has \emph{Property B} with respect to $\{ e_n \}$ if
    for every Laurent polynomial $h$, every contractive $X \in \Mul(\B)$, every $\epsilon > 0$ and $y \in \B$,
    there exists an $M \geq 1$ where for all $\sum b_j$ and $\sum b'_j$ in 
    $\Gamma$ for which 
    $$b_j = b'_j \makebox{  for all  } j \leq M,$$
    we have that 
    $$h\left(\sum b_j X \right)y \approx_{\epsilon}
    h\left(\sum b'_j X \right)y.$$\\

  \item We say that $\Gamma$ has \emph{Property C} with respect to 
    $\{ e_n \}$ if for every Laurent polynomial $h$,
    every contractive $X \in \Mul(\B)$, every $\epsilon > 0$ and
    every $K \geq 1$, there exists an $L \geq 1$ such that for every      
    $\sum b_j$ in $\Gamma$, 
    $$\left \| e_K h\left( \sum b_j X \right) (1 - e_L) \right \| < \epsilon.$$\\ 
  \end{enumerate} 
\end{definition}

\begin{remark}
  Note that the Properties A, B and C described in \Cref{def:Properties-A-B-C} appear asymmetric, but actually
  the opposite versions are immediate consequences stemming from the fact that we may interchange the $z, z^{*}$
  variables in a Laurent polynomial and obtain another Laurent polynomial (i.e., Laurent polynomials are closed
  under adjoints).
\end{remark}

The last 
asymptotic property is of a different nature from the previous properties,
so we single it out in a separate Definition and give it a special name.
Firstly, as preparation, we need to introduce a metalogical concept. 

\begin{df}
  Let $\Gamma_0 \subseteq \Gamma$.
  A statement $\mathfrak{P}(\epsilon)$ is a \emph{finite segment proposition 
    that is sufficient for the subcollection $\Gamma_0$} 
  if $\mathfrak{P}(\epsilon)$ is a proposition
  with free variable $\epsilon$ ranging over $(0, \infty)$ (and possibly other
  free variables) 
  such that if $\sum b_j \in \Gamma$ 
  then $\sum b_j$ belongs to $\Gamma_0$ if there exists a subsequence $\{ M_k \}$ 
  of the positive integers  
  such that
  for all $k$, there exists an $\epsilon > 0$ such that
  $(b_{M_k}, ..., b_{M_{k+1} - 1})$ satisfies $\mathfrak{P}(\epsilon)$, and
  for all $\epsilon > 0$, there exists a $K \geq 1$ such 
  that for all $k \geq K$,
  the finite segment of terms $(b_{M_k}, b_{M_k +1}, ..., b_{M_{k+1}-1})$ satisfies
  the proposition $\mathfrak{P}(\epsilon)$.  
\end{df}

\begin{remark}
  Suppose $\mathfrak{P}(\epsilon)$ is a finite segment proposition such that for
  $\sum b_j \in \Gamma$,
  the series $\sum b_j$ belongs to $\Gamma_0$ if there exists a subsequence $\{ M_k \}$ 
  of the positive integers  
  such that there exists a sequence $\{ \epsilon_k \}$
  in $(0,1)$ with $\sum \epsilon_k < \infty$  where for all $k$, the finite
  segment of terms $(b_{M_k}, b_{M_k +1}, ...,
  b_{M_{k+1} - 1}$ satisfies the proposition $\mathfrak{P}(\epsilon_k)$ (other
  free variables being instantiated in the reasonable way).

  In this case, we say that $\mathfrak{P}(\epsilon)$ \emph{strongly suffices} for $\Gamma_0$.
  We note that so long as the finite segment proposition satisfies the condition that whenever
  some sequence satisfies $\mathfrak{P}(\epsilon)$ it also satisfies $\mathfrak{P}(\epsilon')$ for
  any $\epsilon < \epsilon'$, then if $\mathfrak{P}(\epsilon)$ strongly suffices for $\Gamma_0$, then
  it also suffices for $\Gamma_0$.
\end{remark}

\begin{ex}
Since the above
metalogical concept deviates from the methodology of much of this paper,
we illustrate this with a concrete example.  E.g., 
fix another approximate unit $\{ f_j \}$ for $\B$ for which $f_{j+1} f_j = f_j$
for all $j$.  Let $\Gamma$ be the collection of all series $\sum \alpha_j (f_j - f_{j-1})$ ($f_0 =_{df} 0$) where $\{ \alpha_j \}$ is a bounded sequence in the closed
unit ball $\overline{B(0,1)}$ of the complex plane.
Let $\Gamma_0 \subset \Gamma$ consist of those $\sum \alpha_j (f_j - f_{j-1})$
for which $\{ \alpha_j \}$ is a unit oscillating sequence in $S^1$ with the close 
neighbors property. Then a finite segment proposition $\mathfrak{P}(\epsilon)$
which suffices for $\Gamma_0$ could be the following:
Fix $\epsilon > 0$. 
A finite
segment $(b_{M_k}, ..., b_{M_{k+1}-1}) = (\alpha_{M_k} (f_{M_k} - f_{M_{k-1}}), ...,
\alpha_{M_{k+1}-1} (f_{M_{k+1}-1}- f_{M_{k+1}-2}))$ satisfies the proposition
$\mathfrak{P}(\epsilon)$  
if
\begin{enumerate}
\item $(\alpha_{M_k}, ..., \alpha_{M_{k+1}-1})$ is a finite sequence in $S^1$ that
satisfies conditions (a)-(e) of Definition \ref{df:UnitOscillatingSequence},and 
\item for all $M_k \leq j \leq M_{k+1}$, $|\alpha_{j+1} - \alpha_j | < \epsilon$.
\end{enumerate} 
It is not hard to see that $\sum b_j = 
\sum \alpha_j (f_j - f_{j-1}) \in \Gamma_0$ if and only
if there exists a subsequence $\{ M_k \}$ of the positive integers  
such each segment $(b_{M_k}, ..., b_{M_{k+1}-1})$ satisfies some
 $\mathfrak{P}(\epsilon)$, and for every $\epsilon > 0$, there is a $K \geq 1$,
such that all the finite segments coming after $K$ satisfy $\mathfrak{P}(\epsilon)$.\\ 
\label{ex:MetalogicalEx}
\end{ex}

\begin{definition}
  Let $\mathfrak{P}(\gamma)$ be a finite segment proposition which suffices for
  $\Gamma_0 \subset \Gamma$, where $\gamma$ is a free
  variable ranging over $(0, \infty)$.
The subcollection
 $\Gamma_0$ is said to have \emph{Property X} with respect to $\{ e_n \}$ 
if there exists an $R > 0$ for which    the following statement is true: 

Let $\epsilon > 0$,  $a \in \B_+ - \{ 0 \}$ be contractive, and 
$V \in \Mul(\B)$ be a contractive element for which $\pi(V) \in \C(\B)$ is
unitary.
Let $h_1, h_2, h_3 : S^1 \rightarrow [0,1]$ be nonzero continuous functions and 
$\delta_1 > 0$ such that
$$h_1 h_2 = h_2, \makebox{ } \overline{\osupp(h_3)} \subset  
\osupp((h_2 - \delta_1)_+).$$

There exists $\delta_2 > 0$ such that if $h$ is a Laurent polynomial
for which 
$$| h(\lambda) - h_1(\lambda) | < \delta_2 \makebox{  for all  }
\lambda \in S^1$$
then the following holds:   

For every $L, M \geq 1$, there exist  $L  < L'$,  
$M < M'$, $b_M, b_{M+1}, ..., b_{M'} \in \B$ and $x \in \overline{a \B (e_{L'}
- e_{L})}$ with $\| x \| \leq R$ such that  
for every $\sum b'_j \in \Gamma$ for which
$$b'_j = b_j \makebox{  for all  } M \leq j \leq M',$$
we have that 
$$x h\left( \sum b'_j V \right) x^* \approx_{\epsilon} a.$$ 

Moreover, for any $\gamma > 0$, we may choose $(b_M, b_{M+1}, ..., b_{M'})$
so that it satisfies the proposition $\mathfrak{P}(\gamma)$.
We also require that there is at least one such $\sum b'_j \in \Gamma$ 
for which the above is true.

\end{definition}

Before stating and proving the main theorem of this subsection, 
Theorem \ref{thm:Fullness}, we require
two more definitions.  

\begin{df}
$\Gamma$ is said to be \emph{complete} if a formal series $\sum b_j$ is in $\Gamma$ (and hence converges strictly in $\Mul(\B)$) if and only if for all $k$,
there exists a strictly converging series $\sum b'_j \in \Gamma$ such that
$b_k = b'_k$. 
\end{df}

Also, recall that if $\D$ is a C*-algebra and $x \in \D$, then $x$ is 
\emph{strongly full} in $\D$ if for all $y \in C^*(x) - \{ 0 \}$,
$y$ is full in $\D$, i.e., $Ideal(y) = \D$.

For the convenience of the reader, for 
the main theorem of this subsection, we state explicitly all the hypotheses
that were said, earlier on,
 to be standing assumptions for this section.

\begin{theorem} \label{thm:Fullness} 
Let $\B$ be a $\sigma$-unital, simple, stable C*-algebra.  
Let $\Gamma$ be any collection of strictly converging series of the form
$\sum_{j=1}^{\infty} b_j$, where $b_j \in \B$ for all $j$. Let $\Gamma_0 \subseteq
\Gamma$ be a subcollection such that the finite
segment proposition $\mathfrak{P}(\epsilon)$ strongly suffices for $\Gamma_0$, where
$\epsilon$ is a free variable ranging over $(0, \infty)$.  
Let $\{ e_j \}$ be an approximate unit for $\B$ such that $e_{j+1} e_j = e_j$ 
for all $j$.

Suppose that $\Gamma$ is complete and has Property C with respect to $\{ e_j \}$,
and suppose that $\pi(\Gamma_0)$ consists of unitaries in $\C(\B)$ and $\Gamma_0$
has Property X with respect to $\{ e_j \}$.

Let $V \in \Mul(\B)$ be a contractive element such that $\pi(V)$ is a unitary. 
Then there exists a $\sum b_j \in \Gamma_0$ such that 
$\pi(\sum b_j V)$ is a strongly full unitary in $\C(\B)$.
\end{theorem}

\begin{proof}

For every $k \geq 1$,
let $h_{k, 1, j}, h_{k, 2, j}, h_{k, 3, j} : S^1 \rightarrow [0,1]$
be nonzero 
continuous functions  (for $1 \leq j \leq k$),
and
$\delta_k > 0$ be such that
$$\bigcup_{j=1}^k \osupp(h_{k, 3, j}) = S^1,$$
$$\overline{\osupp(h_{k, 3, j})} \subset \osupp((h_{k, 2, j} - \delta_k)_+),$$
$$h_{k, 1, j} h_{k, 2, j} = h_{k, 2, j},$$
for all $1 \leq j \leq k$, and
$$\max_{1 \leq j \leq k} \diam(\osupp(h_{k, 1, j})) \rightarrow 0$$
as $k \rightarrow \infty$.

Let $\{ h_l \}_{l=1}^{\infty}$ be a sequence of continuous functions
from $S^1$ to $[0,1]$ such that for all $l$, there exists $k, j$
such that $h_l = h_{k, 1, j}$
and for all $k, j$, $h_{k, 1,j}$ occurs infinitely many times as a term
in the sequence $\{ h_l \}_{l=1}^{\infty}$.

Since $\B$ is $\sigma$-unital and stable, we can find a sequence of
norm one 
positive elements $\{ a_n \}_{n=1}^{\infty}$ in $\B$ such that 
$\{ a_n \}$ is an \emph{infinite repeat} of $a_1$, i.e., 
\begin{enumerate}
\item[i.] for all $m, n$, there exists a $y \in \B$ with $y^*y = a_m$
and $yy^* = a_n$, 
\item[ii.] for all $m \neq n$, $a_m \perp a_n$, and
\item[iii.] if $\{ a'_n \}$ is a bounded sequence in $\B$ such that
$a'_n \in \overline{a_n \B a_n}$ for all $n$, then 
$\sum a'_n$ converges strictly in $\Mul(\B)$.
\end{enumerate}
(For example, since $\B$ is stable, viewing $\B$ as $\B \otimes \K$,
let $a_0 \in \B$ be any norm one positive element, and let 
$a_n =_{df} a_0 \otimes e_{n,n}$ for all $n$, where $\{ e_{m,n} \}$
is the standard system of matrix units for $\K$.)

Note that since $\B$ is $\sigma$-unital and simple,
if $\{ a_{n_k} \}$ is any subsequence of $\{ a_n \}$, then
$\sum_{k=1}^{\infty} a_{n_k}$ is a full element of $\Mul(\B)$ (i.e.,
$Ideal(\sum a_{n_k}) = \Mul(\B)$).  In fact, we have that  
there exists a $y' \in \Mul(\B)$ for which 
$y'\left(\sum a_{n_k} \right) (y')^* = 1_{\Mul(\B)}$.  
Moreover, by varying the choice of (contractive positive) $a_1$, we 
may assume that there is always
such a $y'$ where $\| y' \| = 1$.  
We denote the above statement by ``(*)".

Let $\{ \epsilon_k \}$ be a decreasing sequence in $(0,1)$ and
$\{ \epsilon_{k,l} \}$ a (decreasing in $k+l$)
bi-infinite sequence in $(0,1)$ such that
$$\sum_{k=1}^{\infty} \epsilon_k < \infty$$
and
$$\sum_{1 \leq k, l < \infty} \epsilon_{k,l} < \infty.$$
For simplicity, let us assume that $\epsilon_{k,l} = \epsilon_{l,k}$ for all $k, l$.

Before proceeding with Claim 1, we introduce some notation.
Recall that the $n$-th triangular number is given by $\Delta(n) =_{df} \frac{n(n+1)}{2}$.
Suppose that $k$ is an integer for which $\Delta(i) < k < \Delta(i + 1)$.
Notice that $k - \Delta(i) < \Delta(i + 1) - \Delta(i) = i + 1$, and therefore
\begin{equation}
  \label{eq:delta-inequality}
  \Delta(k - \Delta(i)) \le \Delta(i) < k.
\end{equation}

\noindent \emph{Claim 1:}
Since $\Gamma$ is complete, and by repeatedly using that $\Gamma$  has Property C and 
$\Gamma_0$
has Property X, with respect to $\{ e_n \}$, 
we claim that we can find 
a sequence $\{ x_k \}$ of pairwise orthogonal 
elements of $\B$, a sequence $\{ b_l \}$ of elements of $\B$, and a sequence
$\{ \widehat{h}_k \}$ of Laurent polynomials such that the
following statements hold:

\begin{enumerate}
\item $\| x_k \| \leq R$ for all $k$, where $R$ is the constant
from the Definition of Property X. 
\item If $\{ x_{l_j} \}$ is any subsequence of $\{ x_l \}$, 
then  $\sum_{j=1}^{\infty} x_{l_j}$ converges strictly in $\Mul(\B)$.
\item For all $k \geq 1$, if $\Delta(i) < k < \Delta(i+1)$, then
$\widehat{h}_k = \widehat{h}_{\Delta(k-\Delta(i))}$. (See \eqref{eq:delta-inequality}.)
\item \label{item:claim1iv} For all $k \ge 1$, if $k = \Delta(i)$ then 
$$\max_{\lambda \in S^1} |\widehat{h}_k(\lambda) - h_i(\lambda) | < 
\epsilon_k.$$    
\item $\sum b_l$ is a strictly converging series which belongs
to $\Gamma_0$.
\item \label{Oct2020188AM}
 For all $k \geq 1$, there exists a subsequence $\{ x_{l_j} \}$
of $\{ x_l  \}$ such that
\begin{enumerate}
\item $\left\| x_{l_j} \widehat{h}_k\left(\sum_{n=1}^{\infty} b_n V\right) x_{l_s}^* \right\| < \epsilon_{l_j, l_s}$ for all
$j \neq s$, 
\item if $k = \Delta(i)$ for some $i$, then $x_{l_j} \widehat{h}_k\left(\sum_{n=1}^{\infty} b_n V\right) x_{l_j}^* 
\approx_{\epsilon_k} a_{l_j}$  for all $j$, and  
\item if $\Delta(i) < k < \Delta(i + 1)$ for some $i$, then for all $j$,
  \begin{equation*}
    x_{l_j} \widehat{h}_k\left(\sum_{n=1}^{\infty} b_n V\right) x_{l_j}^*
  \approx_{\epsilon_{\Delta(k - \Delta(i))}} a_{l_j}.    
  \end{equation*}
  \end{enumerate}
\end{enumerate}

\noindent Proof of the Claim 1:\\
We will additionally construct three subsequences $\{ L_k \}$, $\{ L'_k \}$ and
$\{ M_k \}$ of the positive integers and a sequence $\{ \delta_{k,2} \}$ of
positive real numbers.  We will always require $0 < \delta_{k,2} < \epsilon_k$
for all $k$.  
The construction is by induction on $k$.\\

\noindent Basis step $k = 1 = \Delta(1)$.
We choose $L_1 = M_1 = 1$.
Since $h_1$ has the form $h_{l_a, 1, j_a}$ for some $l_a$ and $j_a$, 
plug $\epsilon_1$ (for $\epsilon$), $a_1$ (for $a$), $V$,
$h_{l_a,1,j_a}$, $h_{l_a,2,j_a}$, $h_{l_a, 3,j_a}$, $\delta_1$  into 
Property X to 
get a constant $\delta_{1,2} > 0$ (for $\delta_2$). 
We may assume that $\delta_{1,2} < \epsilon_1$.
Since the class of Laurent polynomials is (uniform) norm dense in $C(S^1)$,
choose a Laurent polynomial $\widehat{h}_1$ such that
$$|h_1(\lambda) - \widehat{h}_1(\lambda) | < \delta_{1,2} < \epsilon_1
\makebox{  for all  } \lambda \in S^1.$$
Since $\Gamma_0$ has Property X with respect to $\{ e_n \}$, 
we can choose $M_2 > M_1$, $L'_1 > L_1$, elements $b_1, ..., b_{M_2 - 1} \in \B$,    
and $x_1 \in \overline{a_1 \B (e_{L'_1}  - e_{L_1})}$ with $\norm{x_1} \le R$ such 
that for every $\sum b'_j \in \Gamma$ for which  
$b'_j = b_j$ for all $M_1 \leq j \leq M_2 - 1$, we have that 
$$x_1 \left( \widehat{h}_1\left( \sum b'_j V \right) \right) x_1^*
\approx_{\epsilon_1} a_1.$$
We require that there is at least one such $\sum b'_j \in \Gamma$ as above.
 
We also require that $(b_1, ..., b_{M_2 - 1})$ satisfy the proposition 
$\mathfrak{P}(\epsilon_1)$.\\

\noindent Induction step: Let us assume that 
$\delta_{l,2}$ ($1 \leq l \leq k$),  $\widehat{h}_l$ ($1 \leq l \leq k$), 
$M_{k+1}$, $L'_k > L_k$, $b_1, ..., b_{M_{k+1} - 1}$, and $x_l$ ($1 \leq l
\leq k$)    
have already been chosen.  Recall that we always require
$\delta_{l, 2} < \epsilon_l$ for all $l$.\\  

\noindent \emph{Case 1:} Suppose that $k + 1 = \Delta(i)$ for some $i$. 
Recall that $h_i$ has the form $h_{l_b, 1,j_b}$ for some $s_b, j_b$.
So plug $\epsilon_{k+1}$ (for $\epsilon$), $a_{k+1}$ (for $a$), $V$,
$h_{l_b, 1, j_b}$, $h_{l_b, 2, j_b}$, $h_{l_b, 3, j_b}$, $\delta_{k+1}$ (for
$\delta_1$) into Property X to get a constant $\delta_{k+1, 2} > 0$
(for $\delta_2$).   We may assume that $\delta_{k+1, 2} < \epsilon_{k+1}$.
Since the class of Laurent polynomials is (uniform) norm dense in $C(S^1)$,
choose a Laurent polynomial $\widehat{h}_{k+1}$ such that  
$$| h_i(\lambda) - \widehat{h}_{k+1}(\lambda) | < \delta_{k+1, 2} < 
\epsilon_{k+1} \makebox{  for all  } \lambda \in S^1.$$

Now since $\Gamma$ has Property C with respect to $\{ e_n \}$,
choose  $L_{k+1} > L'_k + 1$ such that 
for all $1 \leq l  \leq k$ and $1 \leq s \leq k+1$,
for every $\sum b'_j \in \Gamma$,
$$\left\| x_l \widehat{h}_{s}\left( \sum b'_j V \right) (1 - e_{L_{k+1} - 1})
\right\|, 
\left\| (1 - e_{L_{k+1} - 1}) \widehat{h}_{s}\left( \sum b'_j V \right) 
x_l^*  \right\| < \epsilon_{k+1, l}.$$

Since $\Gamma_0$ has Property X with respect to $\{ e_n \}$, we can choose 
$M_{k+2} > M_{k+1}$, $L'_{k+1} > L_{k+1}$, elements
$b_{M_{k+1}}, ..., b_{M_{k+2} - 1} \in \B$, and $x_{k+1} \in \overline{a_{k+1} 
\B(e_{L'_{k+1}} - e_{L_{k+1}})}$ with $\norm{x_{k+1}} \le R$
such that  
for every $\sum b'_j \in \Gamma$ for which
$$b'_j = b_j \makebox{ for all } M_{k+1} \leq j \leq M_{k+2} - 1,$$
we have that 
$$x_{k+1} \widehat{h}_{k+1} \left( \sum b'_j V \right) x_{k+1}^*
\approx_{\epsilon_{k+1}} a_{k+1}.$$ 
We require that there exists at least one such $\sum b'_j \in \Gamma$ 
as above, including one where $b'_j = b_j$ for $1 \leq j \leq M_{k+1} - 1$
(from the induction hypothesis, and since $\Gamma$ is complete). 

We also require that 
$(b_{M_{k+1}}, ..., b_{M_{k+2} -1})$ satisfies 
the proposition $\mathfrak{P}(\epsilon_{k+1})$.

Note that for all $1 \leq l \leq k$ and
 $1 \leq s \leq k+1$, for every  
$\sum b'_j \in \Gamma$,
$$\| x_l \widehat{h}_s\left( \sum b'_j V \right) x_{k+1}^* \|, 
\| x_{k+1} \widehat{h}_s \left(\sum b'_j V \right) x_l^* \|
< \epsilon_{k+1, l}.$$\\ 

\noindent \emph{Case 2:} Suppose that $\Delta(i) < k + 1 < \Delta(i + 1)$ for some $i$.  The proof is the same
as in Case 1, but we take $\widehat{h}_{k+1} =_{df} \widehat{h}_{\Delta(k + 1 - \Delta(i))}$
(note that $\widehat{h}_{\Delta(k + 1 - \Delta(i))}$ was chosen in a previous step --- in particular step $\Delta(k + 1 - \Delta(i))$ which is prior to the current step $k + 1$ by \eqref{eq:delta-inequality}),
and there are some minor variations.
For the benefit of the reader, let us nonetheless go through the proof,
proceeding with the aforementioned choice.
For convenience, we will let $S(k+1)$ denote $\Delta(k + 1 - \Delta(i))$.

Since $\Gamma$ has Property C with respect to $\{ e_n \}$, choose 
$L_{k+1} > L'_k + 1$  
such that for all $1 \leq l \leq k$ and $1 \leq s \leq k+1$, 
for all $\sum b'_j \in \Gamma$,
$$\| x_l \widehat{h}_s\left( \sum b'_j V \right) (1 - e_{L_{k+1} -1}) \|,
\| (1 - e_{L_{k+1} -1}) \widehat{h}_s \left(\sum b'_j V \right) x_l^* \|
< \epsilon_{k+1, l}.$$      

In the earlier $S(k+1)$ step of the construction,  
we noted that  $h_{S(k+1)}$ has the form $h_{l_c, 1, j_c}$ for some
$l_c, j_c$, and we plugged in $\epsilon_{S(k+1)}$ (for $\epsilon$), 
$a_{S(k+1)}$ (for $a$), $V$, $h_{l_c, 1, j_c}$, $h_{l_c, 2, j_c}$,
$h_{l_c, 3, j_c}$, $\delta_{S(k+1)}$ (for $\delta_1$) into Property
X to get a constant $\delta_{S(k+1), 2} > 0$ (for $\delta_2$).
We also required that $\delta_{S(k+1), 2} < \epsilon_{S(k+1)}$.
We then chose $\widehat{h}_{S(k+1)}$ so that 
$\max_{\lambda \in S^1} |\widehat{h}_{S(k+1)}(\lambda) - 
h_{S(k+1)}(\lambda) | < \delta_{S(k+1), 2} < \epsilon_{S(k+1)}$. 

So since $\Gamma_0$ has Property X with respect to $\{ e_n \}$, 
and since there exists a $y \in \B$ for which
$y^* y = a_{S(k+1)}$ and $y y^* = a_{k+1}$, 
we can choose
$M_{k+2} > M_{k+1}$, $L'_{k+1} > L_{k+1}$, elements
$b_{M_{k+1}}, ..., b_{M_{k+2} - 1} \in \B$, and 
$x_{k+1} \in \overline{a_{k+1}
\B(e_{L'_{k+1}} - e_{L_{k+1}})}$ with $\norm{x_{k+1}} \le R$
such that
for every $\sum b'_j \in \Gamma$ for which
$$b'_j = b_j \makebox{ for all } M_{k+1} \leq j \leq M_{k+2} - 1,$$
we have that
$$x_{k+1} \widehat{h}_{k+1} \left( \sum b'_j V \right) x_{k+1}^*
\approx_{\epsilon_{S(k+1)}} a_{k+1}.$$
We require that there exists at least one such $\sum b'_j \in \Gamma$
as above, including one where $b'_j = b_j$ for $1 \leq j \leq M_{k+1} - 1$
(from the induction hypothesis, and since $\Gamma$ is complete).

We also require that
$(b_{M_{k+1}}, ..., b_{M_{k+2} -1})$ satisfies
the proposition $\mathfrak{P}(\epsilon_{k+1})$.\\

\noindent This completes the inductive construction.
Finally, we note that $\sum b_l \in \Gamma_0$ since, by our inductive construction, for each $k \ge 1$, $(b_{M_k}, b_{M_{k} + 1}, \ldots , b_{M_{k+1} -1})$ satisfies the finite segment proposition $\mathfrak{P}(\epsilon_k)$ and $\sum \epsilon_k < \infty$.
This proves Claim 1.
\\

Let $U \in \C(\B)$ be the unitary given by $U = \pi(V)$.\\

\noindent \emph{Claim 2:}  For all $m \geq 1$,
$h_m(\pi(\sum_{n=1}^{\infty} b_n) U)$ is a full element of $\C(\B)$.\\

\noindent Proof of Claim 2:
Let $m \geq 1$ be given.
We will now show that $h_m(\pi(\sum_{n=1}^{\infty} b_n) U)$ is full
in $\C(\B)$.
Let $\epsilon > 0$. We may assume that $\epsilon < 1$.
Since each term of the sequence $\{h_l\}_{l=1}^{\infty}$ is repeated infinitely many times, there is some $K \geq 1$ for which $h_K = h_m$ and  
\begin{equation} \label{May2320221AM}\sum_{k=K}^{\infty} \epsilon_k, \makebox{ } 2 \sum_{l=1}^{\infty} \sum_{k =K}^{\infty} \epsilon_{k,l} < \frac{\epsilon}{10(R + 1)^2}. \end{equation}

Let $J =_{df} \Delta(K)$.  Note that $J \geq K$, 
So by Claim 1 \ref{item:claim1iv}, 
\begin{equation} \label{May2320222AM} \max_{\lambda \in S^1} |\widehat{h}_J(\lambda) - h_m(\lambda) |
< \epsilon_J \leq \epsilon_K. \end{equation} (Recall that $h_m = h_K$.)

Also, by Claim 1 \ref{Oct2020188AM}, there exist 
 a subsequence $\{ x_{l_j} \}$, of $\{ x_l \}$, such that 

\begin{equation} \label{May2320223AM}
\| x_{l_j} \widehat{h}_J\left(\sum b_n V \right) x^*_{l_s} \| <
\epsilon_{l_j, l_s} \makebox{ for all  } j \neq s 
\end{equation}
 and
\begin{equation} \label{May2320224AM}
x_{l_j} \widehat{h}_J \left(\sum b_n V \right) x^*_{l_j}   
\approx_{\epsilon_J} a_{l_j} \makebox{  for all  } j. \end{equation}  

We may assume that $l_1 \geq K$.  

Recall that by Claim 1,  
$X = \sum_{j=1}^{\infty} x_{l_j}$ converges strictly to an element of
$\Mul(\B)$ with norm at most $R$. 

Also, recall that the $\{ x_l \}$ are one-sided pairwise orthogonal ($x_l^{*} x_s = 0 = x_l x_s^{*}$ for $s \neq l$) and also 
 $\{ a_l \}$ is pairwise orthogonal.  Hence, 
\begin{eqnarray*}
& & X \widehat{h}_J \left(\sum b_n V \right) X^* \\
& \approx_{\frac{\epsilon}{10(R + 1)^2}}& \sum_{j=1}^{\infty} x_{l_j} \widehat{h}_J
\left(\sum b_n V \right) x_{l_j}^* 
\makebox{   by } (\ref{May2320221AM}), (\ref{May2320223AM})\\
& \approx_{\frac{\epsilon}{10(R+1)^2}} & 
\sum a_{l_j} \makebox{  by  } (\ref{May2320221AM}), (\ref{May2320224AM}).
\end{eqnarray*}
 
By (*), there exists a contractive $y' \in \Mul(\B)$ for which
$y' \left(\sum a_{l_j} \right) (y')^* = 1_{\Mul(\B)}$.
Hence,
$$y'X\widehat{h}_J\left(\sum b_n V\right) X^*(y')^* 
\approx_{\frac{\epsilon}{5(R+1)^2}}  1_{\Mul(\B)}.$$ 

From this,  (\ref{May2320222AM}) and (\ref{May2320221AM}),
\begin{align*}
\pi(y' X) h_m\left(\pi\left(\sum b_n\right) U\right) \pi(X^* (y')^*)
  &\approx_{\frac{\epsilon}{10}} \pi(y' X) \widehat{h}_J\left(\pi\left(\sum b_n \right)
U\right) 
    \pi(X^* (y')^*) \\
  &\approx_{\frac{\epsilon}{5(R+1)^2}} 1_{\C(\B)}.
\end{align*}
                                         
Since $\epsilon < 1$, $\pi(y' X) h_m\left(\pi\left(\sum b_n\right) U\right) \pi(X^* (y')^*)$ is 
invertible in $\C(\B)$, and hence, $h_m\left(\pi\left(\sum b_n\right) U\right)$ is a full element
of $\C(\B)$.

Since $m \geq 1$ was arbitrary, we have proven that for all $m \geq 1$,
$h_m\left(\pi\left(\sum b_n\right) U\right)$ is a full element of $\C(\B)$.
     
\noindent This ends the proof of Claim 2.\\

We now use Claim 2 to complete the proof that 
$\pi(\sum b_n) U$ is a strongly full unitary in $\C(\B)$.  
It suffices to prove that if $f \in C(S^1)$ is any nonzero
positive function, then $Ideal(f(\pi(\sum b_n )U)) = \C(\B)$.

So let $f$ be as above.  By the definition of the sequence $\{ h_l \}$,
there exists $m \geq 1$ and $r > 0$ such that
$$0 \leq rh_m \leq f$$
as functions on $S^1$.  Hence, by the continuous functional calculus,
$$0 \leq r h_m\left(\pi\left(\sum b_n\right)U\right) \leq f\left(\pi\left(\sum 
b_n\right) U\right).$$   
But by Claim 1, $h_m(\pi(\sum b_n) U)$ is a full element of $\C(\B)$.
Hence $f(\pi(\sum b_n) U)$ is a full element of $\C(\B)$.
Since $f$ was arbitrary, $\pi(\sum b_n)U$ is strongly full in
$\C(\B)$, and we are done. 

\end{proof}

\subsection{Proving the Properties}
\label{sec:ProvingTheProperties}

In this subsection, we define, for certain simple C*-algebras, particular 
classes $\Gamma$ and $\Gamma_0 \subseteq \Gamma$, and then show that $\Gamma$ 
possesses Properties A, B and C, and $\Gamma_0$ possesses
Property X, with respect to a relevant approximate
unit.  This, together with our previous results, will lead to $K_1$ injectivity
of certain Paschke dual algebras and the associated uniqueness theorems.

\begin{lemma}
Let $\B$ be a $\sigma$-unital stable C*-algebra, and let $\{ e_n \}$ be an 
approximate unit for $\B$ such that $e_{n+1} e_n = e_n$ for all $n$. 

Suppose that $A, A', A'' \in \C(\B)_+$ are contractive elements and $\delta > 0$
such that
$$A A' = A'$$
and
$$A'' \in her((A' - \delta)_+).$$

Let $A_0 \in \Mul(\B)$ be any contractive lift of $A$, and let $\epsilon > 0$ be
given.

Then for every $L \geq 0$, there exists an $A''_0 \in
\overline{(1 -e_L)\Mul(\B) (1 - e_L)}$ which is a contractive positive
lift of $A''$ such that for all $l \geq 1$,
$$A_0 (A''_0)^{1/l} \approx_{\epsilon} (A''_0)^{1/l} \approx_{\epsilon} (A''_0)^{1/l} A_0.$$
\label{lem:Oct1920186AM}
\end{lemma}

\begin{proof}
The proof is exactly the same as that of \cite[Lemma~4.1]{LN-2020-IEOT}, 
except that
$\{ e_n \}$ need not consist of projections, and for every $n$,
$e_n^{\perp}$ is replaced with $1_{\Mul(\B)} - e_n$, and
$e_n^{\perp} \Mul(\B) e_n^{\perp}$ is replaced with
$\overline{(1_{\Mul(\B)} - e_n) \Mul(\B) (1_{\Mul(\B)} - e_n)}$.
\end{proof}

\vspace*{3ex}

\textbf{We next fix some notation which will be used for the rest of this subsection (\ref{sec:ProvingTheProperties}).  We present this notation here, since it is only for this subsection.} \\

Let $\B$ be a $\sigma$-unital simple stable C*-algebra. 
We will later impose more restrictions on $\B$, but 
do not do so for now.

Let
$\{ e_k \}$ be an approximate unit for $\B$ 
such that $$e_{k+1} e_k = e_k \makebox{ for all  } k.$$
We let $e_0 =_{df} 0$. 

We will later put more restrictions onto $\{ e_k \}$.  In particular,
we will later use Lemma \ref{lem:QuasicentralAU} to require that
$\{ e_k \}$ quasicentralizes a given countable set of elements of 
$\Mul(\B)$.  But for now, we will require no such restrictions.

For all $k$,
let
\begin{equation*}
  r_k \defeq e_k - e_{k-1}.   
\end{equation*}
and for all
$m \leq n$, let
\begin{equation*}
  r_{m,n} \defeq \sum_{k=m}^n r_k = \sum_{k=m}^n (e_k - e_{k-1}) =
  e_n - e_{m-1}.  
\end{equation*}

Note that 
neither the $e_k$ nor the $r_k$  need to be projections.

Recall, also, that we let $\overline{B(0,1)}$ denote the closed unit ball of the
complex plane, i.e., $\overline{B(0,1)} \defeq \{ \alpha \in
\complex : |\alpha| \leq 1 \}$.

Let   
$$\Gamma =_{df} \left\{ \sum_{k=1}^{\infty} \alpha_k r_k :   
\alpha_k \in \overline{B(0,1)}
\makebox{  for all  }k \right\},$$
and let $\Gamma_0 \subseteq \Gamma$ be defined by 
$$\Gamma_0 =_{df} \left\{ \sum_{k=1}^{\infty} \alpha_k r_k : 
 \parbox{40ex}{$\{ \alpha_k \}$ is a unit oscillating sequence in $S^1$
with the close neighbors property} \right\}.$$
Note that in the definitions of $\Gamma$ and $\Gamma_0$ above,
the sums converge strictly in $\Mul(\B)$.\\    

Finally, the finite
segment proposition $\mathfrak{P}(\epsilon)$ (with free variable 
$\epsilon$ ranging over $(0, \infty)$) which is given in  
Example \ref{ex:MetalogicalEx} is sufficient for $\Gamma_0$, as defined above.\\\\

\begin{lemma}
\label{lem:16u}
\label{lem:17b}
\label{lem:PropertiesABC}

$\Gamma$ is complete and 
has Properties A, B and C with respect to $\{ e_n \}$.
\end{lemma}

\begin{proof}
It is not hard to see that $\Gamma$ is complete.

The proofs of Properties A, B and C are exactly the same as those of \cite[Lemmas~4.2 and 4.3]{LN-2020-IEOT}, except
that we do not assume that the $e_n$ are projections and we replace every
occurrence of $e_m^{\perp}$ with $1_{\Mul(\B)} - e_m$ for all $m$.

In fact, the proof is a pretty straightforward operator computation using
the strict topology.
\end{proof}

Towards proving that $\Gamma_0$ has Property X with respect to $\{ e_n \}$,
we firstly prove the following lemma, which is very similar to \cite[Lemma~4.4]{LN-2020-IEOT}, except that we replace
all projections with positive elements and we
have simplified the statement.

\begin{lemma}

Let $h_1, h_{2}, h_{3}: S^1 \rightarrow [0,1]$ be continuous functions and
let $\delta_1 > 0$ be such that
$$h_1 h_{2} = h_{2}$$
and
$$\overline{\osupp(h_{3})} \subset \osupp((h_{2} - \delta_1)_+).$$

Let $\delta_2 > 0$ and $\widehat{h}$ be a Laurent polynomial such that
$$|\widehat{h}(\lambda) - h_1(\lambda) | < \frac{\delta_2}{10}$$
for all $\lambda \in S^1$.

Let $V \in \Mul(\B)$ be a contractive element for which $U =_{df} \pi(V)$
is a unitary in $\C(\B)$.

Then for every $L, M  \geq 1$, 
there exists contractive $A \in (1 -e_{L}) \Mul(\B)_+ (1 - e_L)$ 
which is a lift of $h_{3}(U)$ such that
for every contractive $a \in (\overline{A \B A})_+$,  
there exist $L' > L$, $M' > M$ and
contractive 
$x \in  \overline{r_{L, L'}\B r_{L, L'}}$ for which
$$x \widehat{h} \left( \sum_{j=1}^{\infty} \alpha_j r_j V \right) x^*
\approx_{\delta_2} a$$
for every sequence $\{ \alpha_j \}$  in $\overline{B(0,1)}$ (closed unit
ball of the complex plane) such
that $\alpha_j = 1$ for all $M \leq j \leq M'$.
\label{lem:Oct2020186AM}
\end{lemma}

\begin{proof}

Let $A_0 \in \Mul(\B)_+$ be a contractive lift of $h_1(U)$.

Since $U$ is unitary and because of the conditions on $\hat h$, we know $\hat h(U) \approx_{\frac{\delta_2}{10}} h_1(U)$.
Moreover, since $\hat h$ is a Laurent polynomial, $\hat h(U) = \hat h(\pi(V)) = \pi(\hat h(V))$ and also $h_1(U) = \pi(A_0)$.
Using these facts and since $\Gamma$ has Property A with respect
to $\{ e_n \}$ (by Lemma \ref{lem:16u}), we can
choose $L_1 > L+2$ so that

\begin{equation}
(1 - e_{L_1}) \widehat{h}(V) (1 - e_{L_1}) \approx_{\frac{\delta_2}{10}}
(1 - e_{L_1}) A_0(1 - e_{L_1}) \label{eq:hathV-approx-A0}
\end{equation}

and

\begin{equation}
\widehat{h}\left(\sum_{j=1}^{\infty} \alpha_j r_j V \right)
(1 - e_{L_1-1}) \approx_{\frac{\delta_2}{10}}
\widehat{h}\left(\sum_{j=1}^{\infty} \alpha'_j r_j V \right) (1- e_{L_1-1}) 
\label{equ:Dec320193PM}
\end{equation}
for all sequences $\{ \alpha_j \}$ and $\{ \alpha'_j \}$
in $\overline{B(0,1)}$ such that $\alpha_j = \alpha'_j$ for all $j \geq M$.

By Lemma \ref{lem:Oct1920186AM} (where $A, A', A'', \delta, A_0, \epsilon, L$ are instantiated with $h_1(U), h_2(U),$ $h_3(U), \delta_1, (1 - e_{L_1}) A_0 (1 - e_{L_1}), \frac{\delta_2}{10}, L_1$), 
there exists
an $A \in \overline{(1 - e_{L_1})\Mul(\B)(1 - e_{L_1})}$ which is a 
contractive positive lift of
$h_{3}(U)$ for which
$$(1 - e_{L_1}) A_0 (1 - e_{L_1}) A^{1/l} \approx_{\frac{\delta_2}{10}}
A^{1/l}$$
for all $l \geq 1$.

Hence, if we let $a \in \overline{A \B A}$ be an arbitrary contractive positive
element, then because the previous display holds for all $l \ge 1$,
$$(1 -e_{L_1}) A_0 (1 -e_{L_1}) a^{1/2}
\approx_{\frac{\delta_2}{10}} a^{1/2}.$$
Chaining this with \eqref{eq:hathV-approx-A0} yields
$$(1-e_{L_1}) \widehat{h}(V) (1 -e_{L_1}) a^{1/2}
\approx_{\frac{\delta_2}{5}} a^{1/2}.$$
Therefore,
if we let $y =_{df} a^{1/2}$ then
$$y (1 - e_{L_1}) \widehat{h}(V) (1 - e_{L_1}) y^*
\approx_{\frac{\delta_2}{5}} a.$$

Since $y \in \overline{(1 - e_{L_1}) \B (1 - e_{L_1})}$
and $L_1 > L+2$, 
we can choose
$L' > L_1$ such that if
we define
$$x =_{df} r_{L, L'}y(1 - e_{L_1}) r_{L, L'}$$
then
$$x \widehat{h}(V) x^* \approx_{\frac{\delta_2}{5}} a.$$

Hence, by (\ref{equ:Dec320193PM}),
$$x \widehat{h}\left(\sum_{j=1}^{\infty} \alpha_j r_j V \right) x^*
\approx_{\frac{\delta_2}{2}} a$$
for every sequence $\{ \alpha_j \}$ in $\overline{B(0,1)}$ for which
$\alpha_j = 1$ for all $M \leq j$.
(Note that $(1 - e_{L_1}) r_{L,L'} = (1 - e_{L_1})(e_{L'} - e_{L-1}) = e_{L'} - e_{L_1}$,
and hence $x = x (1 - e_{L_1 - 1})$.)

Note that since $a$ is contractive, $y$ is contractive; and hence,
$x$ is contractive.

Hence,
since $\Gamma$ has Property B with respect to $\{ e_n \}$  
(by
Lemma \ref{lem:16u}),  
we can choose $M' > M$
such that
$$x \widehat{h}\left(\sum_{j=1}^{\infty} \alpha_j r_j V \right) x^*
\approx_{\delta_2}
a$$
for every sequence $\{ \alpha_j \}$ in $\overline{B(0,1)}$ for which
$\alpha_j = 1$ for all $M \leq j \leq M'$.
\end{proof}

The next two short lemmas should be well known, though, for the
convenience of the reader, we indicate the main ideas of the argument.

\begin{lemma}
\label{lem:Kotimesp} 
Let $\D$ be a stable C*-algebra and let $p \in \K$ be a minimal projection.
Let $\iota : \D \rightarrow \D \otimes \K : d \mapsto d \otimes p$.

Then there exists a *-isomorphism $\Phi : \D \rightarrow \D \otimes \K$
such that $\Phi$ and $\iota$ are approximately unitarily equivalent, 
where the unitaries are in $\Mul(\D \otimes \K)$.     
\end{lemma}

\begin{proof}[Idea of the proof]
It is a good exercise to prove that the
 map $\K \rightarrow \K \otimes \K : x \mapsto x \otimes p$ 
is unitarily equivalent to any *-isomorphism $\K \to \K \otimes \K$,
where the unitaries are in $\Mul(\K \otimes \K)$. 
\end{proof}

\begin{lemma} \label{lem:BoundedSubequivalence}
Let $\D$ be a C*-algebra, $a, b, c \in \D$ positive contractive elements, 
$x \in \D$ and $\epsilon > 0$ be such that 
$ab = b$ and $\| x b x^* - c \| < \epsilon$.

Then we can find a $y \in \D$ with $\| y \| \leq \sqrt{1 + \epsilon}$ so that
$\| y a y^* - c \| < \epsilon$ 
\end{lemma}

\begin{proof} Let $y \in \D$ be given by 
$y =_{df} x b^{1/2}$.   So $\| y \|^2 = \| x b x^* \| \leq \| c \| + \epsilon
\leq 1 + \epsilon$. 

Also,  since $y a y^* = x b^{1/2} a b^{1/2} x^* = x b x^*$,
$\| y a y^* - c \| < \epsilon$.
\end{proof}

Recall that if $\D$ is a $\sigma$-unital simple
C*-algebra, then $\D$ has
\emph{strict comparison of positive elements} if for all 
$a, b \in (\D \otimes \K)_+$, 
whenever $d_{\tau}(a) < d_{\tau}(b)$ or $d_{\tau}(b) = \infty$ for all
$\tau \in T(\D)$, we have that $a \preceq b$ (i.e., there exists a
sequence $\{ x_n \}$ in $\D \otimes \K$ for which $x_n b x_n^* \rightarrow b$ 
in norm).

Finally, recall that if $\D$ is a $\sigma$-unital nonunital C*-algebra and
$\tau \in T(\D)$, then there exists an ideal 
$\J_{\tau} \subseteq \Mul(\B)$
which is the norm closure of  $\{ X \in \Mul(\B) : \tau(XX^*) < \infty \}$.
Good references for the study of such ideals are \cite{KNZ-2017-CJM, KNZ-2019-JOT}
and the references contained therein.

\begin{lemma}
  \label{lem:prop_X_aux}
Suppose that, in addition, $\B$ has strict comparison of positive elements
and $T(\B)$ has  finitely many extreme points.

Let $a \in \B$ be a contractive positive element
 and let $\epsilon > 0$ be given.

Let $h_1, h_2, h_3 : S^1 \rightarrow [0,1]$ be continuous functions,
$\delta_1 > 0$ and
$\lambda_1, ..., \lambda_m \in S^1$
such that
$$h_1 h_2 = h_2$$
$$\overline{\osupp(h_3)} \subset \osupp((h_2 - \delta_1)_+)$$
and
the function
$$\lambda \mapsto \sum_{j=1}^m h_3(\lambda_j \lambda)$$
is a full element in $C(S^1)$.

There exists $\delta_2 > 0$ such that if
$\widehat{h}$ is a Laurent polynomial for which
$$|\widehat{h}(\lambda) - h_1(\lambda)| < \frac{\delta_2}{10}$$
for all $\lambda \in S^1$ then the following holds:

Let $V \in \Mul(\B)$ be a contractive element for which
$U = \pi(V)$ is a unitary in $\C(\B)$.  

For every $L, M \geq 1$,  there exist
$M < N_1 < N_2 < ... < N_m$,  $L < L'$,
and
$x \in \overline{a\B r_{L, L'}}$ with $\| x \| \leq 2$ such that
$$x \widehat{h}\left(\sum_{j=1}^{\infty} \alpha_j r_j V\right) x^*
\approx_{\epsilon} a$$
where $\{ \alpha_j \}$ is any sequence in $\overline{B(0,1)}$ such that
$\alpha_j = \lambda_k$ for all $N_{k-1} \leq j < N_k$ and all $1 \leq k \leq m$.
(Here, $N_0 =_{df} M$.)
\label{lem:Oct2020187AM}
\end{lemma}

\begin{proof}
Since $\B$ is stable, we may work with $\B \otimes \K$ in place of
$\B$.

Recall that $\partial T(B \otimes \K)$ is our notation for the extreme boundary 
of $T(\B \otimes \K)$.
For $1 \leq j \leq m$, let $\F_j =_{df} \{ \tau \in \partial T(\B \otimes \K) : 
h_3(\lambda_j U) \notin \pi(\J_{\tau}) \}$, which, by hypothesis, is
a finite set. 
Since the map $\lambda \mapsto \sum_{j=1}^m h_3(\lambda_j \lambda)$ is
norm full in $C(S^1)$, 
$\sum_{j=1}^m h_3(\lambda_j U)$ is norm full in $\C(\B \otimes \K)$. 
So $\bigcup_{j=1}^m \F_j = \partial T(\B \otimes \K)$.

Let $\{ \epsilon_{k,l} \}_{1 \leq k, l \leq m}$ be numbers in $(0,1)$ for
which
\begin{equation} \label{equ:ep_kl} 
\sum_{1 \leq k, l \leq m} \epsilon_{k,l} 
< \frac{\epsilon}{10}.  
\end{equation}

We may assume that $\epsilon_{k,l} = \epsilon_{l,k}$ for all $k, l$.
Replacing $a$ with $(a - \gamma)_+$ for small enough $\gamma > 0$ if
necessary,  we may assume that $d_{\tau}(a) < \infty$ for 
all $\tau \in T(\B)$.\\

\noindent \emph{Step 1:} 
Plug $h_1(\lambda_1 \cdot)$, $h_2(\lambda_1 \cdot )$, $h_3(\lambda_1 \cdot)$,
$\widehat{h}(\lambda_1 \cdot)$ 
(i.e., $h_1, h_2, h_3, \widehat{h}$ all shifted by $\lambda_1$), $\delta_1$, 
$\frac{\epsilon}{100}$ (for $\delta_2$), $V$, $L$ and $M$ into Lemma
\ref{lem:Oct2020186AM} to get a contractive positive operator
 $A_1 \in 
(1 - e_L) \Mul(\B \otimes \K)(1 - e_L)$ such that $\pi(A_1) = h_3(\lambda_1 U)$.
So $A_1 \notin \J_{\tau}$ for all $\tau \in \F_1$, and hence,
$\tau(A_1) = \infty$ for all $\tau \in \F_1$. Let $c \in \overline{A_1 \B A_1}$
be a strictly positive element with norm one.  
Then $d_{\tau}(c) = d_{\tau}(A_1) = \infty$ for all $\tau \in \F_1$.
For all $\tau \in T(\B \otimes \K)$, since the map
$(\B \otimes \K)_+ \mapsto [0, \infty] : d \mapsto d_{\tau}(d)$ is norm lower
semicontinuous, for  $\gamma_1 > 0$ small enough, 
$d_{\tau}((c - \gamma_1)_+) > d_{\tau}(a)$ for all $\tau \in \F_1$.
Let $a_1 =_{df} f_{\gamma_1}(c)$ and $a'_1 =_{df} (c - \gamma_1)_+$.
(See the beginning of \Cref{sec:paschke-dual-simple} for the definition of $f_{\gamma_1}$.)
So
\begin{equation} a_1 \in \overline{A_1 (\B \otimes \K) A_1}, \makebox{ }  
a_1 a'_1 = a'_1 \makebox{ and } d_{\tau}(a'_1) > d_{\tau}(a)
\makebox{   for all }\tau \in \F_1. 
\label{Mar2520226AM}
\end{equation}
By Lemma \ref{lem:Oct2020186AM},
we can find $L'_1 > L$, 
$N_1 > M$ and contractive $x_1  \in  \overline{r_{L'_1, L}(\B \otimes \K)
r_{L'_1, L}}$ for which
$$x_1 \widehat{h}\left( \sum_{j=1}^{\infty} \alpha_j r_j V \right) x_1^*
\approx_{\frac{\epsilon}{10}}  a_1$$
for every sequence $\{ \alpha_j \}$ in $\overline{B(0,1)}$ for which
$\alpha_j = \lambda_1$ for all $M \leq j < N_1$.

By Lemma \ref{lem:Kotimesp}, we can find contractive positive elements
elements $b_1, b'_1 \in \B$ with 
\begin{equation} b_1 b'_1 = b_1 \makebox{ and } d_{\tau}(b'_1 \otimes e_{1,1})
> d_{\tau}(a) \makebox{  for all  } \tau \in \F_1
\label{Mar2520227AM}
\end{equation} 
such that if we replace $x_1$ with 
$y_1 =_{df} (b_1 \otimes e_{1,1})^{1/ n_1}W_1 x_1$, where 
$W_1 \in \Mul(\B \otimes \K)$        
is an appropriate unitary and $n_1 \geq 1$ is a big enough integer, then
$y_1 \in \overline{(b_1 \otimes e_{1,1}) (\B \otimes \K) r_{L'_1, L}}$ is 
contractive and 
\begin{equation}
y_1 \widehat{h}\left( \sum_{j=1}^{\infty} \alpha_j r_j V \right) y_1^*
\approx_{\frac{\epsilon}{10}}  b_1 \otimes e_{1,1}
\end{equation}
for every sequence $\{ \alpha_j \}$ in $\overline{B(0,1)}$ for which
$\alpha_j = \lambda_1$ for all $M \leq j < N_1$.\\

\noindent \emph{Step 2:}  
Since $\Gamma$ has Property C with respect to
$\{ e_n \}$  (by Lemma \ref{lem:17b}),  
choose $L_2 > L'_1 + 3$ 
such that 
\begin{equation} \label{Mar2520228AM}
\| e_{L'_1 + 1} \widehat{h}\left(\sum \alpha_j r_j V \right) 
(1 - e_{L_2 - 2}) \|,  \| (1 - e_{L_2 - 2}) 
\widehat{h} \left( \sum \alpha_j r_j V \right) e_{L'_1 + 1} \|
< \epsilon_{1,2} \end{equation}
for every sequence $\{ \alpha_j \}$ in $\overline{B(0,1)}$.

Plug $h_1(\lambda_2 \cdot)$, $h_2(\lambda_2 \cdot)$, 
$h_3(\lambda_2 \cdot)$, $\widehat{h}(\lambda_2 \cdot)$ (i.e.,
$h_1, h_2, h_3, \widehat{h}$ all shifted by $\lambda_2$), 
$\delta_1$, $\frac{\epsilon}{100}$ (for $\delta_2$), $V$,
$L_2$ (for $L$) and $N_1$ (for $M$) into Lemma \ref{lem:Oct2020186AM}.

Proceeding as in Step 1, we can find $L'_2 > L_2$, $N_2 > N_1$,
contractive positive elements $b_2, b'_2 \in \B \otimes \K$, 
and a contractive element
$y_2 \in \overline{(b_2 \otimes e_{2,2}) (\B \otimes \K) 
r_{L'_2, L_2}}$ such that                  
$$b_2 b'_2 = b'_2 \makebox{  and  } d_{\tau}(b'_2 \otimes e_{2,2}) 
> d_{\tau}(a) 
\makebox{  for all } \tau \in \F_2$$
and
$$y_2 \widehat{h} \left( \sum \alpha_j r_j V \right) y_2^*           
\approx_{\frac{\epsilon}{10}} b_2 \otimes e_{2,2}$$
for every sequence $\{ \alpha_j \}$ in $\overline{B(0,1)}$
for which $\alpha_j = \lambda_2$ for all $N_1 \leq j < N_2$.

Note that by (\ref{Mar2520228AM}),
$$\| y_1  \widehat{h} \left( \sum \alpha_j r_j V \right) y_2^* \|, 
\makebox{  } 
\| y_2 \widehat{h} \left( \sum \alpha_j r_j V \right) y_1^* \| <    
\epsilon_{1,2}$$  
for every  $\{ \alpha_j \}$ in $\overline{B(0,1)}$.\\

\noindent \emph{All Steps:}
By repeating Step 2 multiple times,
we get $M < N_1 < N_2 < ... < N_m$,
$L < L'_1$ and $L'_k + 3 < L_{k+1} < L'_{k+1}$ for all $1 \leq k 
\leq m-1$,  contractive positive elements $b_k, b'_k \in \B$ for all
$1 \leq k \leq m$, and contractive elements $y_k \in 
\overline{(b_k \otimes e_{k,k}) (\B \otimes \K) r_{L'_k, L_k}}$
for all $1 \leq k \leq m$, such that 
\begin{equation}
\label{Mar2520229AM}
b_k b'_k = b'_k \makebox{  and  } d_{\tau}(b'_k \otimes e_{k,k}) > d_{\tau}(a) 
\makebox{  for all } \tau \in \F_k \makebox{ and all } 1 \leq k \leq m, 
\end{equation} 
  
\begin{equation} \label{Mar25202210AM}
y_k \widehat{h} \left( \sum \alpha_j r_j V \right) y_k^*
\approx_{\frac{\epsilon}{10}} b_k \otimes e_{k,k}
\end{equation}
for every sequence $\{ \alpha_j \}$ in $\overline{B(0,1)}$
for which $\alpha_j = \lambda_k$ for all $N_{k-1} \leq j < N_k$
and for all $1 \leq k \leq m$  (we define $N_0 =_{df} M$),
and   
\begin{equation} \label{Mar25202211AM}
\left\| y_k \widehat{h} \left( \sum \alpha_j V \right) y_l^{*} \right\|
< \epsilon_{k,l} \makebox{  for all  } k \neq l 
\makebox{  and  }                             
\end{equation}
all sequences $\{ \alpha_j \}$ in $\overline{B(0,1)}$.

Let $Y =_{df} \sum_{k=1}^m y_k$.  Since the $y_k$s are contractive
and pairwise orthogonal, $\| Y \| \leq 1$.
By (\ref{Mar25202211AM}), (\ref{Mar25202210AM}) and  
(\ref{equ:ep_kl}),
for every sequence $\{ \alpha_j \}$ in $\overline{B(0,1)}$,
for which $\alpha_j = \lambda_k$ for all $N_{k-1} \leq j \leq N_{k}$
and $1 \leq k \leq m$ (recall $N_0 =_{df} M$),
\begin{equation} \label{equ:Mar252022Noon} 
Y \widehat{h}\left(\sum \alpha_j V \right) Y^* \approx_{\frac{\epsilon}{10}} 
\sum_{k=1}^m y_k \widehat{h}\left(\sum \alpha_j V \right) y_k^*
\approx_{\frac{\epsilon}{10}} \sum_{k=1}^m b_k \otimes e_{k,k}. \end{equation}

But by (\ref{Mar2520229AM}) and since $\bigcup_{k=1}^m \F_k = \partial T(\B 
\otimes \K)$, 
$$d_{\tau}\left(\sum_{k=1}^m b'_k \otimes e_{k,k}\right) > d_{\tau}(a)
\makebox{  for all  } \tau \in T(\B \otimes \K).$$ 

Hence, since $\B$ has strict comparison for positive elements,
there exists a sequence $\{ x'_l \}$ in $\B \otimes \K$ such that
$$x'_l \left(\sum_{k=1}^m b'_k \otimes e_{k,k} \right)(x'_l)^*
 \rightarrow a$$
in norm as $l \rightarrow \infty$.
Hence, by  Lemma \ref{lem:BoundedSubequivalence}, 
let $x' \in \B \otimes \K$ with $\|x' \| \leq 2$ be such that 
$$x'\left(\sum_{k=1}^m b_k \otimes e_{k,k}\right) (x')^* 
\approx_{\frac{\epsilon}{10}} a.$$ 
From this and (\ref{equ:Mar252022Noon}),
\begin{equation} \label{equ:Mar2520221PM}
 x'Y\widehat{h}\left(\sum \alpha_j V \right) Y^* (x')^* \approx_{\epsilon}
a \end{equation} 
for every sequence $\{ \alpha_j \}$ in $\overline{B(0,1)}$ for 
which $\alpha_j = \lambda_k$ for all $N_{k-1} \leq j < N_k$ and all
$1 \leq k \leq m$.
      
Note that $\| x'Y \| \leq \| x' \| \| Y \| \leq 2$ and
that equation (\ref{equ:Mar2520221PM}) will still hold if 
we replace $x'Y$ with $a^{1/n} x' Y(r_{L', L})^{1/n}$ for large enough 
$n$ and $L'$. Hence,  (\ref{equ:Mar2520221PM}) still holds with
$x' Y$ replaced with appropriate $z \in \overline{a (\B \otimes \K) r_{L', L}}$ 
with $\| z \| \leq 2$ (and large enough $L'$). 
\end{proof}

\begin{remark}
  By inspection of the proof of \Cref{lem:prop_X_aux}, one may see
  that it is possible to choose $\delta_2$ to depend on $\epsilon$ linearly
  (in the proof, we chose $\frac{\epsilon}{100}$).
  In fact, by varying other constants it may be possible to improve this choice
  further, but here we are not looking for the sharpest estimates.
\end{remark}

\begin{lemma}  $\Gamma_0$ has Property X with respect to  $\{ e_n \}$. 
\label{lem:PropertyX}
\end{lemma}

\begin{proof}
Recall that (as stated in the notation discussion just before
Lemma \ref{lem:17b}) the finite segment
proposition $\mathfrak{P}(\gamma)$ from Example \ref{ex:MetalogicalEx} suffices for $\Gamma_0$, where
 $\gamma$ is a free variable ranging
over $(0, \infty)$. 

That $\Gamma_0$ has Property X basicly follows from
Lemma \ref{lem:Oct2020187AM}.  We just need to make 
sure that, given $\epsilon > 0$, we  modify the
finite segment of Lemma \ref{lem:Oct2020187AM} so that it
satisfies the proposition $\mathfrak{P}(\epsilon)$.
The finite segment (or output) of Lemma \ref{lem:Oct2020187AM} is   
\begin{equation}
\label{Mar2520223PM}
(\lambda_1 r_{N_0}, \lambda_1 r_{N_0 +1}, ..., \lambda_1 r_{N_1 - 1}, 
\lambda_2 r_{N_1}, ..., \lambda_2 r_{N_2 - 1}, ...., \lambda_m r_{N_{m-1}},
... \lambda_m r_{N_{m} - 1}).\end{equation}         
(Recall that $N_0 =_{df} M$.)
Note that in the above finite segment, the coefficent $\lambda_1$ repeats
at least 
$N_1 - N_0$ times, the coefficient $\lambda_2$   
repeats at least $N_2 - N_1$ times, and in general, for $1 \leq j \leq m$,
the coefficient $\lambda_j$ repeats at least $N_j - N_{j-1}$ times.

We can always rearrange or 
add more terms to the sequence $\{ \lambda_k \}_{k=1}^m$
before plugging it into the hypothesis of 
Lemma \ref{lem:Oct2020187AM}. 
In the output of Lemma \ref{lem:Oct2020187AM}, which is
the finite segment (\ref{Mar2520223PM}), 
we can always vary the coefficients by a small amount (and thus remove 
some of the repetitions) 
and still retain the conclusion of Lemma \ref{lem:Oct2020187AM}
(note that a Laurent polynomial is uniformly continuous on the closed unit
ball of any \cstar-algebra). 
Thus, we can ensure that in the conclusion of Lemma \ref{lem:Oct2020187AM},
the finite segement (\ref{Mar2520223PM}) is replaced with
the finite segment
$$(\lambda'_1 r_{N_0}, \lambda'_1 r_{N_0 +1}, ..., \lambda'_H r_{N_M})$$
where $|\lambda'_j - \lambda'_{j+1}| < \epsilon$ for all $1 \leq j \leq  
H - 1$ (thus ensuring the close neighbors property) and
where the sequence $\{ \lambda'_j \}_{j=1}^H$ satisfies 
conditions (a) to (e) of Definition \ref{df:UnitOscillatingSequence} (where
$\{ \lambda'_j \}_{j=1}^H = \{ \alpha_j \}_{j=M_k}^{M_{k+1} - 1}$; thus
ensuring unit oscillation).  Thus, the proposition $\mathfrak{P}(\epsilon)$
is satisfied.  
\end{proof}

\vspace*{5ex}

We are now ready to prove the main theorem of this subsection.
Unlike in the previous results of this subsection, we will explicitly state
all the hypotheses, many of which were standing assumptions for most of
this subsection (as stated in the discussion before Lemma \ref{lem:17b}).

\begin{theorem}
Let $\B$ be a $\sigma$-unital, stable, finite, simple C*-algebra with
strict comparison for positive elements and with $T(\B)$ having finitely many 
extreme points. 
Let $S \subseteq \Mul(\B)$ a countable set,
and let $\{ e_n \}$ be an approximate unit for
$\B$ that  satisfies the hypotheses of Lemma \ref{lem:QuasicentralAU} 
with respect to the given $S$. 

For all $n \geq 1$, let $r_n =_{df} e_n - e_{n-1}$ (with $e_0 =_{df} 0$).
Let $U \in \C(\B)$ be a unitary.

Then there exists a sequence $\{ \alpha_n \}$ in $S^1$ such that
$\pi\left(\sum \alpha_n r_n\right)$ is a unitary in $\pi(S)'$ such that
$\pi\left(\sum \alpha_n r_n\right)$ is homotopy equivalent to $1$ in $\pi(S)'$
(i.e., there is a norm continuous path of unitaries in $\pi(S)'$
that has $\pi\left(\sum \alpha_n r_n\right)$ and $1$ as its endpoints) 
such that $\pi\left(\sum \alpha_n r_n\right)U$ is a strongly full unitary
in $\C(\B)$.
\label{thm:Mar2520225PM}
\end{theorem}
 
\begin{proof}
Let $\Gamma$ be the set of strictly converging series in $\Mul(\B)$
that is given by
$$\Gamma =_{df} \left\{ \sum_{j=1}^{\infty} \beta_j r_j : \beta_j \in
\overline{B(0,1)} \makebox{  for all } j \right\}$$
and let $\Gamma_0 \subseteq \Gamma$ be given by
$$\Gamma_0 =_{df} \left\{ \sum_{j=1}^{\infty} \beta_j r_j :
\parbox{40ex}{$\{ \beta_j \}$ is a unit oscillating sequence in $S^1$
with the close neighbors property} \right\}.$$

By Lemmas \ref{lem:AUUnitary} and \ref{lem:QuasicentralAU}, 
every element of $\pi(\Gamma_0)$ is a unitary in $\pi(S)'$.
By Lemma \ref{lem:HomotopyTrivialInCommutant}, every such unitary
is homotopy trivial in $\pi(S)'$. I.e., every element  $x \in 
\pi(\Gamma_0)$
is a unitary in $\pi(S)'$, and there is a norm continuous  
path of unitaries in $\pi(S)'$ connecting $x$ to $1$.

By Lemma \ref{lem:17b}, $\Gamma$ is complete and
has Property C with respect to
$\{ e_n \}$.  
Recall that the finite segment proposition
$\mathfrak{P}(\gamma)$ given in Example \ref{ex:MetalogicalEx} suffices for $\Gamma_0$.
By Lemma \ref{lem:PropertyX}, $\Gamma_0$ has Property X  
with respect to $\{ e_n \}$.  By Theorem \ref{thm:Fullness},
there exists a series $\sum \alpha_j r_j \in \Gamma_0$  
such that $\pi\left( \sum \alpha_j r_j \right)U$ is a strongly
full unitary in $\C(\B)$.   
\end{proof}

\subsection{$\ko$-injectivity and uniqueness}

Having laid the groundwork, we are finally in a position to prove our main theorem concerning the $\ko$-injectivity of the Paschke dual algebra $\paschkedual{\A}{\B}$, and obtain the interesting uniqueness result.\\

\begin{theorem} \label{thm:K1InjectivityStrictComp} 
  \label{thm:k1-injectivity-simple-nuclear-strict-comparison}
  Suppose that $\A,\B$ are separable simple \cstar-algebras with $\A$ unital 
and nuclear,  with $\B$ stable and having strict comparison of 
positive elements, and
with $T(\B)$ having finitely many extreme points. 

  Then $\paschkedual{\A}{\B}$ is $\ko$-injective.
  Moreover, for each $n \in \nat$, the map
  \begin{equation*}
    U(\mat_n \otimes \paschkedual{\A}{\B})/U(\mat_n \otimes \paschkedual{\A}{\B})_0 \to U(\mat_{2n} \otimes \paschkedual{\A}{\B})/U(\mat_{2n} \otimes \paschkedual{\A}{\B})_0
  \end{equation*}
  given by
  \begin{equation*} [u] \mapsto [u \oplus 1] \end{equation*}
  is injective.
\end{theorem}

\begin{proof}
Since $\B$ has strict comparison for positive elements,
$\B$ has the corona factorization property (see \cite{OPR-2012-IMRN}).

Let $\phi : \A  \rightarrow \Mul(\B)$ be a unital trivial absorbing extension.
Thus, we may realize the Paschke dual algebra $\A^d_{\B}$ as
$\A^d_{\B} =  (\pi \circ \phi(\A))' \in \C(\B)$.
Let $U \in (\pi \circ \phi(\A))' \in \C(\B)$ be an arbitrary unitary.
By Theorem \ref{thm:Mar2520225PM}, there exists a unitary $W \in 
(\pi \circ \phi(\A))'$ such that
$W \sim_h 1$ in $(\pi \circ \phi(\A))'$ and
$W U$ is strongly full in $\C(\B)$.  (In Theorem \ref{thm:Mar2520225PM},
take $S = \phi(S_0)$ where $S_0 \subset \A$ is a countable 
dense subset.)
Since $U$ was arbitrary, the result follows from 
Theorem \ref{thm:k1-injectivity-generic}.
\end{proof}

Now we obtain a uniqueness result, which is
 our primary generalization of the BDF essential 
codimension result.

\maintheorem

\begin{proof}
  Since $\B$ has strict comparison, $\B$ has the corona factorization property (see \cite{OPR-2012-IMRN}).
  Hence, $\phi$ and $\psi$ are both (unital) absorbing extensions (e.g., see the paragraph preceding \Cref{thm:k1-injectivity-generic}).

  By Theorem \ref{thm:K1InjectivityStrictComp},
$\A^d_{\B}$ is $\ko$-injective, so the result follows from Theorem \ref{thm:MainUniqueness}.
\end{proof}

\bibliographystyle{amsalpha}
\bibliography{references.bib}

\newcommand{\etalchar}[1]{$^{#1}$}
\providecommand{\bysame}{\leavevmode\hbox to3em{\hrulefill}\thinspace}
\providecommand{\MR}{\relax\ifhmode\unskip\space\fi MR }
\providecommand{\MRhref}[2]{%
  \href{http://www.ams.org/mathscinet-getitem?mr=#1}{#2}
}
\providecommand{\href}[2]{#2}
\begin{thebibliography}{{Weg}93}

\bibitem[Arv77]{Arv-1977-DMJ}
William Arveson, \emph{Notes on extensions of {{\(C^*\)}}-algebras}, Duke Math.
  J. \textbf{44} (1977), 329--355 (English).

\bibitem[Arv07]{Arv-2007-PNASU}
William Arveson, \emph{Diagonals of normal operators with finite spectrum},
  Proc. Natl. Acad. Sci. USA \textbf{104} (2007), no.~4, 1152--1158 (English).

\bibitem[BCP{\etalchar{+}}06]{BCP+-2006-Agatoeo}
Moulay-Tahar Benameur, Alan~L. Carey, John Phillips, Adam Rennie, Fyodor~A.
  Sukochev, and Krzysztof~P. Wojciechowski, \emph{An analytic approach to
  spectral flow in von {N}eumann algebras}, Analysis, geometry and topology of
  elliptic operators (Matthias Lesch, Bernhelm Boo-Bavnbek, Slawomir Klimek,
  and Weiping Zhang, eds.), World Sci. Publ., Hackensack, NJ, 2006,
  pp.~297--352. \MR{2246773}

\bibitem[BDF73]{BDF-1973-PoaCoOT}
Lawrence~G. Brown, Ronald~George Douglas, and Peter~Arthur Fillmore,
  \emph{Unitary equivalence modulo the compact operators and extensions of
  {$C^{\ast} $}-algebras}, {Proceedings of a Conference on Operator Theory}
  (Peter~Arthur Fillmore, ed.), Lecture Notes in Mathematics, vol. 345,
  Springer, Berlin, 1973, pp.~58--128. \MR{0380478 (52 \#1378)}

\bibitem[BJ15]{BJ-2015-TAMS}
Marcin {Bownik} and John {Jasper}, \emph{{The Schur-Horn theorem for operators
  with finite spectrum.}}, {Trans. Am. Math. Soc.} \textbf{367} (2015), no.~7,
  5099--5140 (English).

\bibitem[BL12]{BL-2012-CJM}
Lawrence~G. Brown and Hyun~Ho Lee, \emph{Homotopy classification of projections
  in the corona algebra of a non-simple {$C^*$}-algebra}, Canad. J. Math.
  \textbf{64} (2012), no.~4, 755--777. \MR{2957229}

\bibitem[Bla98]{Bla-1998}
Bruce Blackadar, \emph{{$K$}-theory for operator algebras}, second ed.,
  Mathematical Sciences Research Institute Publications, vol.~5, Cambridge
  University Press, Cambridge, 1998. \MR{1656031}

\bibitem[BRR08]{BRR-2008-JNG}
Etienne Blanchard, Randi Rohde, and Mikael R\o{r}dam, \emph{Properly infinite
  {$C(X)$}-algebras and {$K_1$}-injectivity}, J. Noncommut. Geom. \textbf{2}
  (2008), no.~3, 263--282. \MR{2411419}

\bibitem[Dav96]{Dav-1996}
Kenneth~R. Davidson, \emph{{$C^*$}-algebras by example}, Fields Institute
  Monographs, vol.~6, American Mathematical Society, Providence, RI, 1996.
  \MR{1402012}

\bibitem[DE01]{DE-2001-KT}
Marius Dadarlat and S\o{r}en Eilers, \emph{Asymptotic unitary equivalence in
  {$KK$}-theory}, $K$-Theory \textbf{23} (2001), no.~4, 305--322. \MR{1860859}

\bibitem[EG96]{EG-1996-AM}
George~A. Elliott and Guihua Gong, \emph{On the classification of
  {{\(C^*\)}}-algebras of real rank zero. {II}}, Ann. Math. (2) \textbf{144}
  (1996), no.~3, 497--610 (English).

\bibitem[EK01]{EK-2001-PJM}
George~A. Elliott and Dan Kucerovsky, \emph{An abstract
  {V}oiculescu-{B}rown-{D}ouglas-{F}illmore absorption theorem}, Pacific J.
  Math. \textbf{198} (2001), no.~2, 385--409. \MR{1835515}

\bibitem[ERS11]{ERS-2011-AJM}
George~A. Elliott, Leonel Robert, and Luis Santiago, \emph{The cone of lower
  semicontinuous traces on a {C}{{\(^{\ast}\)}}-algebra}, Am. J. Math.
  \textbf{133} (2011), no.~4, 969--1005 (English).

\bibitem[{Gab}16]{Gab-2016-PJM}
James {Gabe}, \emph{{A note on nonunital absorbing extensions.}}, {Pac. J.
  Math.} \textbf{284} (2016), no.~2, 383--393 (English).

\bibitem[GN19]{GN-2019-AM}
Thierry {Giordano} and Ping~W. {Ng}, \emph{{A relative bicommutant theorem: the
  stable case of Pedersen's question}}, {Adv. Math.} \textbf{342} (2019), 1--13
  (English).

\bibitem[GR20]{GR-2020-GMJ}
James Gabe and Efren Ruiz, \emph{The unital {Ext}-groups and classification of
  {{\(C^*\)}}-algebras}, Glasg. Math. J. \textbf{62} (2020), no.~1, 201--231
  (English).

\bibitem[Hig87]{Hig-1987-PJM}
Nigel Higson, \emph{A characterization of {$KK$}-theory}, Pacific J. Math.
  \textbf{126} (1987), no.~2, 253--276. \MR{869779}

\bibitem[Jas13]{Jas-2013-JFA}
John Jasper, \emph{{The Schur--Horn theorem for operators with three point
  spectrum}}, J. Funct. Anal. \textbf{265} (2013), no.~8, 1494--1521.
  \MR{3079227}

\bibitem[JT91]{JT-1991}
Kjeld~Knudsen Jensen and Klaus Thomsen, \emph{Elements of {$KK$}-theory},
  Mathematics: Theory \& Applications, Birkh\"auser Boston, Inc., Boston, MA,
  1991. \MR{1124848}

\bibitem[Kad02]{Kad-2002-PNASUa}
Richard~V. Kadison, \emph{{The Pythagorean Theorem II: the infinite discrete
  case}}, Proc. Natl. Acad. Sci. USA \textbf{99} (2002), no.~8, 5217--5222.

\bibitem[Kas80]{Kas-1980-JOT}
G.~G. Kasparov, \emph{Hilbert {$C^{\ast} $}-modules: theorems of {S}tinespring
  and {V}oiculescu}, J. Operator Theory \textbf{4} (1980), no.~1, 133--150.
  \MR{587371}

\bibitem[Kas81]{Kas-1981-MUI}
\bysame, \emph{The operator {K}-functor and extensions of {C}*-algebras}, Math.
  USSR, Izv. \textbf{16} (1981), 513--572 (English).

\bibitem[KL17]{KL-2017-IEOT}
Victor Kaftal and Jireh Loreaux, \emph{Kadison's pythagorean theorem and
  essential codimension}, Integr. Equ. Oper. Theory \textbf{87} (2017),
  565--580.

\bibitem[KN06]{KN-2006-HJM}
Dan Kucerovsky and P.~W. Ng, \emph{The corona factorization property and
  approximate unitary equivalence}, Houston J. Math. \textbf{32} (2006), no.~2,
  531--550. \MR{2219330}

\bibitem[KNZ09]{KNZ-2009-JFA}
Victor Kaftal, Ping~Wong Ng, and Shuang Zhang, \emph{Strong sums of projections
  in von {N}eumann factors}, J. Funct. Anal. \textbf{257} (2009), no.~8,
  2497--2529. \MR{2555011 (2011b:46096)}

\bibitem[KNZ17]{KNZ-2017-CJM}
Victor {Kaftal}, Ping~Wong {Ng}, and Shuang {Zhang}, \emph{{Strict comparison
  of positive elements in multiplier algebras}}, {Can. J. Math.} \textbf{69}
  (2017), no.~2, 373--407 (English).

\bibitem[KNZ19]{KNZ-2019-JOT}
Victor Kaftal, P.~W. Ng, and Shuang Zhang, \emph{Purely infinite corona
  algebras}, J. Oper. Theory \textbf{82} (2019), no.~2, 307--355 (English).

\bibitem[Lee11]{Lee-2011-JFA}
Hyun~Ho Lee, \emph{Proper asymptotic unitary equivalence in {$KK$}-theory and
  projection lifting from the corona algebra}, J. Funct. Anal. \textbf{260}
  (2011), no.~1, 135--145. \MR{2733573}

\bibitem[Lee13]{Lee-2013-JFA}
\bysame, \emph{Deformation of a projection in the multiplier algebra and
  projection lifting from the corona algebra of a non-simple {$C^*$}-algebra},
  J. Funct. Anal. \textbf{265} (2013), no.~6, 926--940. \MR{3067791}

\bibitem[Lee18]{Lee-2018-JMAA}
\bysame, \emph{Homotopy classification of homogeneous projections in the corona
  algebra of {$C(X,B)$} for a graph {$X$}}, J. Math. Anal. Appl. \textbf{461}
  (2018), no.~2, 1404--1415. \MR{3765498}

\bibitem[Lin91]{Lin-1991-PAMS}
Huaxin Lin, \emph{Simple {{\(C^*\)}}-algebras with continuous scales and simple
  corona algebras}, Proc. Am. Math. Soc. \textbf{112} (1991), no.~3, 871--880
  (English).

\bibitem[Lin01]{Lin-2001-CJM}
\bysame, \emph{Classification of simple tracially {AF} {{\(C^*\)}}-algebras},
  Can. J. Math. \textbf{53} (2001), no.~1, 161--194 (English).

\bibitem[Lin02]{Lin-2002-JOT}
\bysame, \emph{Stable approximate unitary equivalence of homomorphisms}, J.
  Operator Theory \textbf{47} (2002), no.~2, 343--378. \MR{1911851}

\bibitem[{Lin}05]{Lin-2005-CJM}
Huaxin {Lin}, \emph{{Extension by simple \(C^*\)-algebras: quasidiagonal
  extensions}}, {Can. J. Math.} \textbf{57} (2005), no.~2, 351--399 (English).

\bibitem[LN20]{LN-2020-IEOT}
Jireh Loreaux and P.~W. Ng, \emph{Remarks on essential codimension}, Integral
  Equations Oper. Theory \textbf{92} (2020), no.~4, 1--35.

\bibitem[Lor19]{Lor-2019-JOT}
Jireh Loreaux, \emph{Restricted diagonalization of finite spectrum normal
  operators and a theorem of arveson}, Journal of Operator Theory \textbf{81}
  (2019), no.~2, 257--272.

\bibitem[Ng18]{Ng-2018-NYJM}
Ping~W. Ng, \emph{A double commutant theorem for the corona algebra of a
  {R}azak algebra}, New York J. Math. \textbf{24} (2018), 157--165.
  \MR{3761942}

\bibitem[OPR11]{OPR-2011-TAMS}
Eduard {Ortega}, Francesc {Perera}, and Mikael {R{\o}rdam}, \emph{{The Corona
  factorization property and refinement monoids}}, {Trans. Am. Math. Soc.}
  \textbf{363} (2011), no.~9, 4505--4525 (English).

\bibitem[OPR12]{OPR-2012-IMRN}
\bysame, \emph{{The corona factorization property, stability, and the Cuntz
  semigroup of a \(C^*\)-algebra}}, {Int. Math. Res. Not.} \textbf{2012}
  (2012), no.~1, 34--66 (English).

\bibitem[Pas81]{Pas-1981-PJM}
William~L. Paschke, \emph{{$K$}-theory for commutants in the {C}alkin algebra},
  Pacific J. Math. \textbf{95} (1981), no.~2, 427--434. \MR{632196}

\bibitem[{Ped}90]{Ped-1990}
Gert~K. {Pedersen}, \emph{{The corona construction}}, {Operator theory, Proc.
  GPOTS-Wabash Conf., Indianapolis/IN 1988, Pitman Res. Notes Math. Ser. 225,
  49-92 (1990).}, 1990.

\bibitem[R{\o}r04]{Ror-2004-IJM}
Mikael R{\o}rdam, \emph{The stable and the real rank of {{\(\mathcal
  Z\)}}-absorbing {{\(C^*\)}}-algebras}, Int. J. Math. \textbf{15} (2004),
  no.~10, 1065--1084 (English).

\bibitem[Tho01]{Tho-2001-PAMS}
Klaus Thomsen, \emph{On absorbing extensions}, Proc. Amer. Math. Soc.
  \textbf{129} (2001), no.~5, 1409--1417. \MR{1814167}

\bibitem[TT15]{TT-2015-CMB}
Aaron~Peter Tikuisis and Andrew Toms, \emph{On the structure of {Cuntz}
  semigroups in (possibly) nonunital {{\(C^\ast\)}}-algebras}, Can. Math. Bull.
  \textbf{58} (2015), no.~2, 402--414 (English).

\bibitem[Val83]{Val-1983-PJM}
Alain Valette, \emph{A remark on the {K}asparov groups {${\rm Ext}(A,\,B)$}},
  Pacific J. Math. \textbf{109} (1983), no.~1, 247--255. \MR{716300}

\bibitem[Voi76]{Voi-1976-RRMPA}
Dan Voiculescu, \emph{A non-commutative {W}eyl-von {N}eumann theorem}, Rev.
  Roumaine Math. Pures Appl. \textbf{21} (1976), no.~1, 97--113. \MR{0415338}

\bibitem[{Weg}93]{Weg-1993}
Niels~Erik {Wegge-Olsen}, \emph{\({K}\)-theory and \({C}^*\)-algebras: a
  friendly approach}, Oxford: Oxford University Press, 1993 (English).

\end{thebibliography}

\end{document}